\documentclass[onefignum,onetabnum,hidelinks]{siamart171218}

\usepackage{lipsum}
\usepackage{amsfonts}
\usepackage{graphicx}
\usepackage{epstopdf}
\usepackage{algorithmic}
\ifpdf
  \DeclareGraphicsExtensions{.eps,.pdf,.png,.jpg}
\else
  \DeclareGraphicsExtensions{.eps}
\fi

\newsiamremark{remark}{Remark}
\newsiamremark{hypothesis}{Hypothesis}
\crefname{hypothesis}{Hypothesis}{Hypotheses}
\newsiamthm{claim}{Claim}

\title{Convergence analysis of a Crank-Nicolson Galerkin method for an inverse source problem for parabolic equations with boundary observations
}

\author{Dinh Nho H\`ao\thanks{Hanoi Institute of Mathematics, VAST, 18 Hoang Quoc Viet Road, Hanoi, Vietnam ({\tt hao@math.ac.vn})}
	\and Tran Nhan Tam Quyen\thanks{Institute for Numerical and Applied Mathematics, University of Goettingen, Lotzestr. 16-18, 37083 Goettingen, Germany ({\tt quyen.tran@uni-goettingen.de}) (Corresponding author)}
	\and Nguyen Thanh Son\thanks{ICTEAM Institute, Universit\'e catholique de Louvain, 1348 Louvain-la-Neuve, Belgium and Department of Mathematics and Informatics, Thai Nguyen University of Sciences, 25000 Thai Nguyen, Vietnam ({\tt thanh.son.nguyen@uclouvain.be})}
}

\usepackage{amsopn}

\ifpdf
\hypersetup{
}
\fi

\def\<{(}
\def\>{)_{\partial\Omega}}

\def\limk{\lim_{k\to\infty}}
\def\limsupk{\limsup_{k\to\infty}}
\def\liminfk{\liminf_{k\to\infty}}

\begin{document}

\maketitle

\begin{abstract}
This work is devoted to an inverse problem of identifying a source term depending on both spatial and time variables in a parabolic equation from single Cauchy data on a part of the boundary. A Crank-Nicolson Galerkin method is applied to the least squares functional with an quadratic stabilizing penalty term. The convergence of finite dimensional regularized approximations to the sought source as measurement noise levels and  mesh sizes approach to zero with an appropriate regularization parameter is proved. Moreover, under a suitable source condition, an error bound and corresponding convergence rates are proved. Finally, several numerical experiments are
presented to illustrate the theoretical findings.
\end{abstract}



\section{Introduction}
The problem of identifying a source in a heat transfer or diffusion process has got attention of many researchers during last years. This problem leads to determining a term in the right hand side of parabolic equations from some observations of the solution which is well known to be ill-posed.  For surveys on the subject, we refer the reader to the books  \cite{cannon1984,haob,isakov1990,isakov2006,prilepko-book2000}, the recent papers \cite{hao2017, prilepko-tkachenko2003c} and the references therein.

Although there have been many papers devoted to the source identification problems with observations in the whole domain or at the final moment,  those with boundary observations are quite few. Furthermore, the sought term depends either on the spatial variable as in  \cite{cannon1968,cannon-duchateau1998,cannon-ewing,cannon-lin,choulli_yamamoto2004,choulli_yamamoto2006,engl_scherzer_yamamoto,hao4,hao5,haob,trong_alain_nam,yamamoto1993,yamamoto1994}, or only on the time variable as  in \cite{hasanov_pektas2014}. In this paper, we consider the problem of determining the right hand side depending on both spatial and time variables by a variational method. Indeed, 
let $\Omega$ be an open bounded connected subset of $\mathbb{R}^d,
d\ge 1$ with boundary $\partial \Omega$ and $T>0$ be a given constant. We investigate the problem of identifying the {\it source term} $f=f(x,t)$ in the Robin boundary value problem for the parabolic equation
\begin{equation}\label{ct1}
\begin{aligned}
&\frac{\partial u}{\partial t} (x,t) + \mathcal{L}u(x,t) =f(x,t) \text{~in~} \Omega_T := \Omega \times (0,T],  \\
&\frac{\partial u(x,t)}{\partial\vec{n}} +\sigma(x,t)u(x,t) = g(x,t) \text{~on~} \mathcal{S} := \partial\Omega \times (0,T],\\
&u(x,0) = q(x) \text{~in~} \Omega
\end{aligned}
\end{equation}
from a partial boundary measurement $z_\delta := z_\delta(x,t) \in L^2(\Sigma)$ of the solution $u(x,t)$ on the surface $\Sigma := \Gamma\times(0,T) \subset\mathcal{S}$ satisfying
\begin{align}\label{15-4-19ct1}
\|Z - z_\delta\|_{L^2(\Sigma)} \le \delta,
\end{align}
where $Z := u_{|\Sigma}$, $\Gamma$ is a relatively open subset of $\partial\Omega$ and the positive constant $\delta$ stands for the measurement error. 

In \cref{ct1} $\mathcal{L}$ is a time-dependent, second order self-adjoint elliptic operator of the form
$$\mathcal{L}u(x,t) := - \sum_{i,j=1}^d \frac{\partial}{\partial x_i}\left( a_{ij}(x,t)\frac{\partial u(x,t)}{\partial x_j} \right) +b(x,t)u(x,t),$$
where $A :=  \left(a_{ij}\right)_{1\le i,j\le d} \in {C(\overline{\Omega_T})}^{d \times d}$ is a symmetric diffusion matrix satisfying the uniformly elliptic condition
\begin{align}\label{15-1-17ct2}
A(x,t)\xi \cdot\xi = \sum_{i,j=1}^d a_{ij}(x,t)\xi_i\xi_j \ge \underline{a} \sum_{i,j=1}^d |\xi_i|^2 \mbox{~in~} \overline{\Omega_T}
\end{align}
for all $\xi = \left(\xi_i\right)_{1\le i\le d} \in \mathbb{R}^d$ with some constant $\underline{a} >0$ and $b(x,t)\in C(\overline{\Omega_T})$ is a non-negative function. The vector $\vec{n} :=\vec{n}(x,t)$ is the unit outward normal on $\mathcal{S}$ and
$$\frac{\partial u(x,t)}{\partial\vec{n}} := A(x,t)\nabla u(x,t) \cdot\vec{n}$$
with $\nabla u(x,t) := \nabla_x u(x,t) = \left( \frac{\partial u(x,t)}{\partial x_1},\ldots,\frac{\partial u(x,t)}{\partial x_d}\right)$. In addition, the functions $q\in H^1(\Omega)$, $g\in C(\overline{\mathcal{S}})$ and $\sigma\in C(\overline{\mathcal{S}})$ with $\sigma(x,t)\ge 0$ in $\overline{\mathcal{S}}$ are assumed to be given. The source term $f=f(x,t)$ is sought in the space $L^2(\Omega_T)$.

The contents of this paper are as follows: For any fixed $f\in L^2(\Omega_T)$ let $u=u(f) \in \mathcal{W}(0,T)$ denote the unique weak solution of the system \cref{ct1}, see \Cref{definition} for the definition of related functional spaces. Adopting the output least squares method combined with the Tikhonov regularization, we consider the (unique) minimizer of the minimization problem
$$
\min_{f\in L^2(\Omega_T)}J_{\rho,\delta}(f) \quad\mbox{with}\quad J_{\rho,\delta}(f) := \|u(f)-z_\delta\|^2_{L^2(\Sigma)} + \rho\|f-f^*\|^2_{L^2(\Omega_T)} \eqno \left(\mathcal{P}_{\rho,\delta}\right)
$$
as a reconstruction, where $\rho\in (0,1)$ is the regularization parameter and $f^*$ is an a priori estimate of the exact source which is identified. We mention that $f^*$ plays the role of a selection criterion. In practice it is not easy to
estimate $f^*$, however an a priori knowledge of the identified source improves the quality of reconstructions (cf.\ \cite{EHN96,kir11}).

For discretization we employ the Crank-Nicolson Galerkin method, where the finite dimensional space $\mathcal{V}_h^1 := \big\{\varphi_h\in C\left(\overline\Omega\right)
~|~{\varphi_h}_{|T} \in \mathcal{P}_1(T) \mbox{~for all~}
T\in \mathcal{T}_h\big\}$ of piecewise linear,
continuous finite elements is used to discretize the state with respect to the spatial variable. Further, to discretize the state with respect to the time variable, we divide the time interval $(0,T)$ into $M$ equal subintervals and introduce a time step $\tau := T/M$ together with time levels
\begin{align}\label{9-12-19ct1}
t^n := n\tau \mbox{~with~} n \in I_0:= \{0,1,\ldots,M\}.
\end{align}
As a result, the state $u(f)$ is then approximated by the finite sequence $(U^n_{h,\tau}(f))_{n=0}^M$ in which for each $n\in I_0$ the element 
$U^n_{h,\tau}(f) \in \mathcal{V}_h^1.$ 
With these notions at hand, we examine the discrete regularized problem corresponding to $\left(\mathcal{P}_{\rho,\delta}\right)$ i.e. the following strictly convex minimization
problem
$$
\min_{f\in L^2(\Omega_T)}J_{\rho,\delta,h,\tau}(f) \eqno \left(\mathcal{P}_{\rho,\delta,h,\tau}\right)
$$
with
$$
 J_{\rho,\delta,h,\tau}(f) := \sum_{n=1}^M \int_{t^{n-1}}^{t^n} \|U^n_{h,\tau}(f)-z_\delta\|^2_{L^2(\Gamma)}dt + \rho\|f-f^*\|^2_{L^2(\Omega_T)}
$$
which admits a unique solution $f_{\rho,\delta,h,\tau}$ obeying the relation (see \Cref{4-6-19ct1})
\begin{align}\label{4-6-19ct3}
{f_{\rho,\delta,h,\tau}}_{|\Omega\times(t^{n-1},t^n]} = f^*- \rho^{-1}P^{n-1}_{h,\tau}(f_{\rho,\delta,h,\tau}) 
\end{align}
for any $n\in I := \{1,\ldots,M\}$, where $(P^n_{h,\tau}(f))_{n=0}^M$ is the approximation of the adjoint state $p(f)$. Using the variational discretization concept introduced in \cite{hin05}, the minimizer automatically belongs to the finite dimensional space 
$$
\mathcal{V}_{h,\tau}^{1,0} := \left\{\Phi\in L^2\left( 0,T;\mathcal{V}^1_h\right)
~|~{\Phi}_{|\Omega\times(t^{n-1},t^n]} := \varphi^n_h \in \mathcal{V}^1_h \mbox{~for all~} n \in I\right\},
$$
thus a discretization
of the admissible set $L^2(\Omega_T)$ can be avoided. 

As $\delta,h,\tau \to 0$ and with an appropriate a priori regularization parameter choice $\rho=\rho(\delta,h,\tau) \to 0$, we show in \Cref{convergence} that
the whole sequence $\big(f_{\rho,\delta,h,\tau}\big)_{\rho>0}$ converges in the $L^2(\Omega_T)$-norm to the unique $f^*$-minimum-norm solution $f^\dag$ of the
identification problem  defined by
$$
f^\dag = \arg \min_{f \in \left\{ f \in L^2(\Omega_T) ~|~ u(f)_{|\Sigma} =  Z \right\}}  \|f-f^*\|_{L^2(\Omega_T)}. \eqno\left(\mathcal{IP}\right)
$$
The corresponding state sequence then converges in the $L^2(0,T; H^1(\Omega))$-norm to the state $u(f^\dag)$.

\Cref{rate} is devoted to convergence rates for the discretized problem, where we first show that if $f\in \left\{ f \in L^2(\Omega_T) ~|~ u(f)_{|\Sigma} =  Z \right\}$ and there exists a function $w \in L^2(\Sigma)$ such that $f = F(w)+f^*$, where $F(w)$ is the unique weak solution of the parabolic equation
\begin{equation}\label{10-12-19ct1}
\begin{aligned}
&-\frac{\partial F}{\partial t} (x,t) + \mathcal{L}F(x,t) =0 \text{~in~} \Omega_T,  \\
&\frac{\partial F(x,t)}{\partial\vec{n}} +\sigma(x,t)F(x,t) = w\chi_\Sigma  \text{~on~} \mathcal{S}, \\
&F(x,T) = 0 \text{~in~} \Omega
\end{aligned}
\end{equation}
with $\chi_\Sigma$ being the characteristic function of $\Sigma\subset\mathcal{S}$, then $f = f^\dag$. Furthermore, if the data appearing in the system \cref{ct1} are regular enough the error bound
\begin{align*}
\|f_{\rho,\delta,h,\tau} - f^\dag\|^2_{L^2(\Omega_T)} \le C\left( h^3\rho^{-1} + \tau^2h^{-1}\rho^{-1} +\delta +\rho + \delta^2\rho^{-1}\right)
\end{align*}
is established, that yields the convergence rate \begin{align*}
\|f_{\rho,\delta,h,\tau} - f^\dag\|_{L^2(\Omega_T)} = \mathcal{O}(\sqrt{\delta})
\end{align*}
with the suitable choice of the parameters $h$, $\tau$ and $\rho$.

For the numerical solution of the discrete regularized problem $\big(\mathcal{P}_{\rho,\delta,h,\tau}\big)$ we in \Cref{NumerEx} utilize a conjugate gradient algorithm. Numerical studies are presented for several cases where the sought source is smooth and discontinuous as well, that illustrates the efficiency of the proposed approach.

To conclude this introduction we wish to mention that to the best of our knowledge, although there have been many papers devoted to source identification problems for parabolic equations, we however have not yet found investigations
on the discretization analysis for those with boundary observations --- which is more realistic from the practical point of view, a fact that motivated the research presented in the paper.

Regarding the identification problem in elliptic equations utilizing {\it boundary} measurements, we here would like to comment briefly some previously published works. In \cite{XiZo05,XuZo15} the authors used finite element methods to numerically recover the fluxes on the inaccessible boundary $\Gamma_i$ from measurement data of the state on the accessible boundary $\Gamma_a$, while the problem of identifying the Robin coefficient on $\Gamma_i$ is also considered in \cite{XuZo15a}. In \cite{HHQ19,HKQ18,HQ19,quyenANM} some source and coefficient identification problems have been investigated, and in \cite{quyen18} the problem of simultaneously identifying the source term and coefficients from {\it distributed} observations. A survey of the optimal control problems for parabolic equations can be found in the book \cite{troltzsh}, where distributed data is also taken into account.

Throughout the paper we use the
standard notion of Sobolev spaces $H^1(\Omega)$, $H^2(\Omega)$, $W^{k,p}(\Omega)$, etc. from, for example, \cite{ada75}. 

\section{Problem setting and preliminaries}\label{definition}

To formulate the identification problem, we first give some notations \cite{wolka}. 
Let $H^1(\Omega)^*$ be the dual space of $H^1(\Omega)$, we use the notation
\begin{align*}
\mathcal{W}(0,T) := \left\{ w\in L^2\left( 0,T; H^1(\Omega)\right) ~\big|~ \frac{\partial w}{\partial t} \in L^2\left( 0,T; H^1(\Omega)^*\right) \right\}.
\end{align*}
It is a Banach space  equipped with the norm
\begin{align}\label{15-1-17ct4}
\|w\|_{\mathcal{W}(0,T)} := \left( \|w\|^2_{L^2\left( 0,T; H^1(\Omega)\right)} + \left\|\frac{\partial w}{\partial t}\right\|^2_{L^2\left( 0,T; H^1(\Omega)^*\right)}\right)^{1/2}.
\end{align}
We note that, since $\mathcal{W}(0,T)$ with respect to the norm \cref{15-1-17ct4} is a closed subspace of the reflexive space $L^2\left( 0,T; H^1(\Omega)\right)\times L^2\left( 0,T; H^1(\Omega)^*\right)$, it is itself reflexive. 

\subsection{Direct and inverse problems}\label{D-I problems}

For considering \Cref{ct1}, we set
\begin{equation*}
\begin{aligned}
&a(v,w) \\
~\quad &:=\int_\Omega  A(x,t)\nabla v(x)\cdot\nabla w(x)dx + \int_\Omega
b(x,t)v(x)w(x)dx +\int_{\partial\Omega}\sigma(x,t)v(x)w(x)dx
\end{aligned}
\end{equation*}
where $t\in (0,T]$ and $v,w \in H^1(\Omega)$. Then, for each $f\in L^2(\Omega_T)$ \Cref{ct1}
defines a unique weak solution $u := u(x,t;f) := u(f)$ in the sense that $u(f) \in \mathcal{W}(0,T)$ with $u(x,0) = q(x)$ for a.e. $x \in \Omega$ and the following variational equation is satisfied (cf.\ \cite{troltzsh,wolka})
\begin{align}\label{17-10-16ct2}
\left\langle \frac{\partial u}{\partial t}, \varphi \right\rangle_{\left( H^1(\Omega)^*, H^1(\Omega)\right) } + a(u,\varphi)  = (f,\varphi)_{L^2(\Omega)} + ( g, \varphi)_{L^2(\partial\Omega)} 
\end{align}
for all $\varphi\in H^1(\Omega)$ and a.e. $t\in (0,T]$. Furthermore, the estimate
\begin{align}\label{17-10-16ct4}
\|u\|_{\mathcal{W}(0,T)} \le C \left( \|f\|_{L^2(\Omega_T)} + \|g\|_{L^2(\mathcal{S})} + \|q\|_{L^2(\Omega)}\right)
\end{align}
holds, where $C$ is a positive constant independent of $f$, $g$ and $q$.
Therefore, we can define the {\it source-to-state} operator
$$u: L^2(\Omega_T) \to \mathcal{W}(0,T)$$
which maps each $f\in L^2(\Omega_T)$ to the unique weak solution $u:= u(f)$ of \Cref{ct1}. The inverse problem is stated as follows:
\begin{quote} \it
\begin{center}
	Given the boundary data $Z = {u}_{|\Sigma}$ of the exact solution $u$, find an element $f\in L^2(\Omega_T)$ such that $u(f)_{|\Sigma} =Z$.
	\end{center}
\end{quote}

\subsection{Variational method}\label{ssec:VariationMeth}
In practice only the observation $z_\delta \in L^2(\Sigma)$ of the exact $Z$
with an error level
\begin{align*}
\left\|Z - z_\delta \right\|_{L^2(\Sigma)} \le \delta, \quad \delta>0
\end{align*}
is available. Hence, our problem is to reconstruct an element $f\in L^2(\Omega_T)$ in \cref{ct1} from noisy observation $z_\delta$ of $Z$. For this purpose we use the standard least squares method with Tikhonov regularization, i.e. we consider a minimizer of the minimization problem
$
\left(\mathcal{P}_{\rho,\delta}\right)
$
as a reconstruction. 

\begin{remark}
	In case $ \inf_{(x,t)\in\Omega_T} b(x,t) >0$ or $\inf_{(x,t)\in\mathcal{S}}\sigma(x,t) >0$ the expression
	$a(v,w)$
	generates a scalar product on the space
	$ H^1(\Omega)$ equivalent to the usual one, i.e. there exist positive constants $C_1, C_2$ such that (cf.\ \cite{Mik78,Pechstein})
	\begin{align}\label{18-10-16ct1}
	C_1\|\varphi\|^2_{H^1(\Omega)} \le a(\varphi,\varphi) \le C_2\|\varphi\|^2_{H^1(\Omega)}
	\end{align}
	for all $\varphi\in H^1(\Omega)$ and $t\in (0,T]$.
	
	Now we assume that $b =\sigma=0$. A change of the variable $u=e^tv$, \Cref{ct1} has the form
	\begin{equation*}
	\begin{aligned}
	&\frac{\partial v}{\partial t} (x,t) - \sum_{i,j=1}^d \frac{\partial}{\partial x_i}\left( a_{ij}(x,t)\frac{\partial v(x,t)}{\partial x_j} \right) + v(x,t) = e^{-t}f(x,t) \text{~in~} \Omega_T,  \\
	&\frac{\partial v(x,t)}{\partial\vec{n}} = e^{-t}g(x,t) \text{~on~} \mathcal{S} \quad\mbox{and}\quad
	v(x,0) = q(x) \text{~in~} \Omega.
	\end{aligned}
	\end{equation*}
	In the sequel we thus consider the case $ \inf_{(x,t)\in\Omega_T} b(x,t) >0$ or $\inf_{(x,t)\in\mathcal{S}}\sigma(x,t) >0$ only. All results in present paper are still valid for the case $b =\sigma
	=0$.
\end{remark}

\begin{remark}\label{weakly conv.}
If the sequence $\left( f_k\right)\subset L^2(\Omega_T)$
		weakly converges in $L^2(\Omega_T)$ to an element $f$, then the sequence $\left( u(f_k)\right)$ weakly converges in $\mathcal{W}(0,T)$ to $u(f)$.
\end{remark}

For each $f\in L^2(\Omega_T)$ we now consider the {\it adjoint} problem
\begin{equation}\label{18-1-17ct1}
\begin{aligned}
&-\frac{\partial p}{\partial t} (x,t) + \mathcal{L}p(x,t) =0 \text{~in~} \Omega_T,  \\
&\frac{\partial p(x,t)}{\partial\vec{n}} +\sigma(x,t)p(x,t) = \left( u(x,t;f)-z_\delta(x,t)\right)\chi_\Sigma  \text{~on~} \mathcal{S}, \\
&p(x,T) = 0 \text{~in~} \Omega,
\end{aligned}
\end{equation}
where $\chi_\Sigma$ is the characteristic function of $\Sigma\subset\mathcal{S}$, i.e. $\chi_\Sigma(x,t) =1$ if $(x,t)\in\Sigma$ and equals to zero otherwise. A function $p \in \mathcal{W}(0,T)$ is said to be a weak solution to this problem, if $p(x,T) = 0$ for a.e.\ $x \in \Omega$ and 
\begin{align}\label{17-10-16ct2-ad}
-\left\langle \frac{\partial p}{\partial t}, \varphi \right\rangle_{\left( H^1(\Omega)^*, H^1(\Omega)\right) } + a(p,\varphi)  =  \int_\Gamma (u(x,t;f) - z_\delta(x,t)) \varphi(x) dx
\end{align} 
for all $\varphi\in H^1(\Omega)$ and a.e. $t\in (0,T]$. Since $u \in \mathcal{W}(0,T)$, the boundary value $\left( u(x,t;f)-z_\delta(x,t)\right)\chi_\Sigma$ belongs to $L^2(\Sigma)$ and by changing the time direction we see that \cref{18-1-17ct1} attains a unique weak solution $p(x,t;f) := p(f) \in \mathcal{W}(0,T)$.

We close this section with the following result.

\begin{theorem}\label{existance}
	The minimization problem $\left(\mathcal{P}_{\rho,\delta}\right)$ attains a unique minimizer $f_{\rho,\delta}$ which satisfies the equation
	\begin{align}\label{19-10-16ct6}
	f_{\rho,\delta}=f^*- \rho^{-1} p(f_{\rho,\delta}).
	\end{align}
\end{theorem}

\begin{proof}
Let $(f_k)$ be a minimizing sequence of the problem $\left(\mathcal{P}_{\rho,\delta}\right)$, i.e.
		\begin{align*}
		\limk J_{\rho,\delta}(f_k) = \inf_{f\in L^2(\Omega_T)} J_{\rho,\delta}(f).
		\end{align*}
		The sequence $(f_k)$ is then bounded in the $L^2(\Omega_T)$-norm. A subsequence not relabeled and an element $f\in L^2(\Omega_T)$ exist such that $(f_k)$ weakly converges to $f$ in $L^2(\Omega_T)$. Due to \Cref{weakly conv.} and the continuity of the mapping $L^2(0,T; H^1(\Omega)) \to L^2(0,T; L^2(\Gamma)) = L^2(\Sigma)$, it follows that $u(f_k)$ weakly converges to $u(f)$ in $L^2(\Sigma)$. And since the $L^2$-norm is weakly lower semi-continuous, is so the functional $J_{\rho,\delta}(\cdot)$. Furthermore, it is clear that $J_{\rho,\delta}(f)$ is strictly convex, hence it attains a unique minimizer. Meanwhile, the equation \cref{19-10-16ct6} follows directly from the first order optimality condition, which finishes the proof.
\end{proof}

\section{Crank-Nicolson Galerkin discretization}\label{4-6-19ct1}

In this section we present the {\it Crank-Nicolson Galerkin} method (see, e.g., \cite{thome}) to discretize the regularized minimization problem in finite dimensional spaces.

Let $\left(\mathcal{T}_h\right)_{0<h<1}$ be a regular family of
quasi-uniform triangulations of the domain $\overline{\Omega}$ with the mesh size $h$ (cf.\ e.g., \cite{brenner}). We first consider the finite element space $\mathcal{V}_h^1$
and recall the Cl\'ement interpolation operator
$\Pi_h: L^1(\Omega) \rightarrow \mathcal{V}_h^1$
which satisfies the following properties
\begin{equation}\label{23/10:ct2}
\lim_{h\to 0} \left\| v - \Pi_h v
\right\|_{H^k(\Omega)} =0 \enskip \mbox{for all} \enskip k \in \{0,
1\}
\end{equation}
and
\begin{equation}\label{23/5:ct1}
\left\| v - \Pi_h v \right\|_{H^k(\Omega)} \le
Ch^{l-k} \| v\|_{H^l(\Omega)}
\end{equation}
for $0 \le k < l \le 2$ (see \cite{Clement,Bernardi2,scott_zhang}). We also mention that for all $v\in \mathcal{V}^1_h$, the inverse inequality
\begin{align}\label{17-10-16ct5**}
\|v\|_{H^1(\Omega)} \le Ch^{-1}\|v\|_{L^2(\Omega)}
\end{align}
holds true. Further,
for all $v\in H^1(\Omega)$ and $T\in\mathcal{T}_h$ there holds the local estimate (\cite{nochetto})
\begin{align}\label{17-10-16ct5***}
\|v\|_{L^2(\partial T)} \le C\left( h^{-1/2}\|v\|_{L^2(T)} + h^{1/2}\|\nabla v\|_{L^2(T)}\right).
\end{align} 

To discretize the state functions with respect to the time variable, we consider the partition \cref{9-12-19ct1} of the interval $(0,T)$. 
For a continuous function $\zeta : [0,T] \to \mathbb{R}$ and $t^n\in [0,T]$, we denote by
$$\zeta^n := \zeta(t^n).$$
Then, we set
\begin{align*}
&a^n(v,w) \\
~\quad &=\int_\Omega A^n(x)\nabla v(x)\cdot\nabla w(x) dx + \int_\Omega
b^n(x)v(x)w(x)dx +\int_{\partial\Omega}\sigma^n(x)v(x)w(x)dx
\end{align*}
for all $v,w\in H^1(\Omega)$. 

Linking to the above partition, we introduce the constant piecewise, discontinuous interpolation operator
\begin{align}\label{1-2-17ct1}
\pi_\tau : L^2(0,T) \to \mathcal{V}_\tau^0 := \left\{\varphi \in L^2(0,T)
~|~{\varphi}_{|(t^{n-1},t^n]} = \mbox{~constant,~} \forall n\in I\right\}
\end{align}
which is defined for each $\eta\in L^2(0,T)$ by
$${\pi_\tau\eta(t)}_{|(t^{n-1},t^n]} := {\eta_\tau(t)}_{|(t^{n-1},t^n]} :=\frac{1}{\tau}\int_{t^{n-1}}^{t^n} \eta(s)ds :=\overline{\eta}^n,\quad\forall n\in I.$$
By Proposition 9 of \cite[pp.\ 129]{royden}, we get the limit
\begin{align}\label{1-2-17ct1*}
\lim_{\tau\to 0}\|\eta-\eta_\tau\|_{L^2(0,T)}=0
\end{align}
and the estimate (cf.\ \cite[Proposition 5.1]{kroner})
\begin{align}\label{1-2-17ct1**}
\|\eta-\eta_\tau\|_{L^2(0,T)} \le C\tau \left\|\frac{d\eta}{dt}\right\|_{L^2(0,T)}
\end{align}
in case $\eta\in H^1(0,T)$.

For a sequence $(w^k) \in L^2(\Omega)$ we respectively introduce  the backward difference quotient and the mean as follows
$$
\partial w^k := \frac{w^k-w^{k-1}}{\tau} \quad\mbox{and}\quad \overline{\partial} w^k := \frac{w^k+w^{k-1}}{2}.
$$

We are in the position to present the Crank-Nicolson Galerkin method applied to \cref{ct1}. Let $q_h := \Pi_h q \in \mathcal{V}_h^1$, find $U^n(f):=U^n_{h,\tau}(f)\in \mathcal{V}^1_h$ such that $U^0(f) = q_h$ and 
\begin{equation}\label{18-1-17ct5*}
\begin{aligned}
\left( \partial U^n(f), \varphi_h \right)_{L^2(\Omega)} + a^n(\overline{\partial}U^n(f),\varphi_h) &= (\overline{f}^n,\varphi_h)_{L^2(\Omega)} + ( g^n, \varphi_h)_{L^2(\partial\Omega)}
\end{aligned}
\end{equation}
for all $\varphi_h\in \mathcal{V}^1_h$, $n\in I$.

\begin{definition}
	The mapping $u_{h,\tau}: L^2(\Omega_T)\to \left( \mathcal{V}_h^1\right)^M$ defined for each $f\in L^2(\Omega_T)$ by
	$$u_{h,\tau}(f) := \left( U^1(f),\ldots,U^M(f)\right) \in \left( \mathcal{V}_h^1\right)^M $$
	is called the {\it discrete source-to-state} operator.
\end{definition}
The operator $u_{h,\tau}$ is Fr\'echet differentiable on $L^2(\Omega_T)$. For each $f\in L^2(\Omega_T)$
in the direction $\xi\in L^2(\Omega_T)$  its differential $u_{h,\tau}'(f)\xi$ is $\left( \widehat{U}^1(\xi),\ldots,\widehat{U}^M(\xi)\right)$, where $\widehat{U}^n(\xi) := \widehat{U}^n_{h,\tau}(\xi) \in \mathcal{V}^1_h$ is defined by the  variational equation
\begin{equation}\label{18-1-17ct5**}
\begin{aligned}
\left( \partial \widehat{U}^n(\xi), \varphi_h \right)_{L^2(\Omega)} + a^n(\overline{\partial}\widehat{U}^n(\xi),\varphi_h) &= (\overline{\xi}^n,\varphi_h)_{L^2(\Omega)}, \quad \forall \varphi_h\in \mathcal{V}^1_h,~ n\in I\\
\widehat{U}^0(\xi) &=0.
\end{aligned}
\end{equation}
The adjoint problem \cref{18-1-17ct1} is discretized via the process that for $n\in I_0$ the element $P^n(f) := P^n_{h,\tau}(f) \in \mathcal{V}^1_h$ satisfies the system
\begin{equation}\label{18-1-17ct5***}
\begin{aligned}
-\tau \left( \partial P^n(f), \varphi_h \right)_{L^2(\Omega)} + \tau a^n(\overline{\partial}P^n(f),\varphi_h) = \int_{t^{n-1}}^{t^n}( U^n(f)-z_\delta, \varphi_h)_{L^2(\Gamma)}dt
\end{aligned}
\end{equation}
for all $\varphi_h\in \mathcal{V}^1_h$ with
$P^M(f) =0$. Note that from $P^M(f)=0$ we can compute $P^{M-1}(f)$ due to \cref{18-1-17ct5***}, and so $P^{M-2}(f)$, ..., $P^1(f),P^0(f)$.

To begin, we recall the discrete Gronwall inequality.

\begin{lemma}
	Assume that $(Y_n)_n$,~ $(X_n)_n$,~ and $(\alpha_n)_n$ are non-negative sequences such that
	$$Y_n \le X_n + \sum_{m=0}^{n-1} \alpha_m Y_m \quad\mbox{for all}\quad n\ge 0.$$
	Then
	$$Y_n \le X_n + \sum_{m=0}^{n-1} \alpha_m X_m e^{\sum_{s=m+1}^{n-1} \alpha_s}$$
	for all $n\ge 0$.
\end{lemma}

\begin{lemma}\label{disc.solu}
	
	(i) There holds the estimate
	\begin{equation}\label{19-1-17ct4}
	\max_{n\in I_0}\left\| U^n(f) \right\|^2_{H^1(\Omega)}\le C\left( \| q \|^2_{H^1(\Omega)}+ \|g\|^2_{L^\infty\left( 0,T; L^2(\partial\Omega)\right) }  + \|f\|^2_{L^2(\Omega_T)} \right).
	\end{equation}
	
	(ii) The inequalities
	\begin{align}\label{19-1-17ct4*}
	& \sum_{n=1}^M \int_{t^{n-1}}^{t^n} \left|\left( \partial U^n(f),\theta\right)_{L^2(\Omega)}\right|dt \notag\\
	&~\quad \le C\left( \| q \|_{H^1(\Omega)}+ \|g\|_{L^\infty\left( 0,T; L^2(\partial\Omega)\right) }  + \|f\|_{L^2(\Omega_T)} \right)\|\theta\|_{L^2\left(0,T; H^1(\Omega)\right) }
	\end{align}
	for all $\theta\in L^2\left(0,T; H^1(\Omega) \right)$ and
	\begin{align}\label{19-1-17ct4***}
	&\sum_{n=1}^M \tau \sup_{\|\varphi\|_{H^1(\Omega)}\le 1}\left|\left( \partial U^n(f),\varphi\right)_{L^2(\Omega)}\right|^2 \notag\\
	&~\quad \le C\left( \| q \|^2_{H^1(\Omega)}+ \|g\|^2_{L^\infty\left( 0,T; L^2(\partial\Omega)\right) }  + \|f\|^2_{L^2(\Omega_T)} \right).
	\end{align}
	hold true.
	
	(iii) The limit
	\begin{equation}\label{19-1-17ct4**}
	\lim_{h,\tau\to 0} \sum_{n=1}^M \int_{t^{n-1}}^{t^n} \left( U^n(f)-U^{n-1}(f), \theta\right)_{H^1(\Omega)}dt =0
	\end{equation}
	is satisfied for all $\theta\in L^2\left(0,T; H^1(\Omega) \right)$.
\end{lemma}

\begin{proof}
	(i) Taking $\varphi_h = 2\tau^2 \partial U^n(f)$ in \cref{18-1-17ct5*}, we have
	\begin{align*}
	&2\|U^n(f)-U^{n-1}(f)\|^2_{L^2(\Omega)} +\tau a^n(U^n(f),U^n(f)) \\
	&~\quad = \tau a^n(U^{n-1}(f),U^{n-1}(f)) + 2\tau(\overline{f}^n,U^n(f)-U^{n-1}(f))_{L^2(\Omega)} \\
	&~\qquad + 2\tau(  g^n, U^n(f)-U^{n-1}(f))_{L^2(\partial\Omega)}.
	\end{align*}
	By \cref{18-10-16ct1}, we thus get
	\begin{align*}
	&2\|U^n(f)-U^{n-1}(f)\|^2_{L^2(\Omega)} + C_1\tau \|U^n(f)\|^2_{H^1(\Omega)}\\
	&~\quad \le C_2\tau \|U^{n-1}(f)\|^2_{H^1(\Omega)}  + 2\tau\|\overline{f}^n\|_{L^2(\Omega)}\|U^n(f)-U^{n-1}(f)\|_{L^2(\Omega)} \\
	&~\qquad + 2\tau\| g^n\|_{L^2(\partial\Omega)} \|U^n(f)-U^{n-1}(f)\|_{L^2(\partial\Omega)}.
	\end{align*}
	For an arbitrary $\epsilon>0$, an application of Young's inequality yields that
	\begin{align*}
	&2\tau\|\overline{f}^n\|_{L^2(\Omega)}\|U^n(f)-U^{n-1}(f)\|_{L^2(\Omega)}\\
	&~\quad \le \epsilon\|U^n(f)-U^{n-1}(f)\|^2_{L^2(\Omega)}+\frac{\tau^2}{\epsilon}\|\overline{f}^n\|^2_{L^2(\Omega)} \\
	&~\quad \le \epsilon\|U^n(f)-U^{n-1}(f)\|^2_{L^2(\Omega)}+\frac{\tau}{\epsilon}\int_{t^{n-1}}^{t^n}\|f\|^2_{L^2(\Omega)} dt.
	\end{align*}
	Meanwhile, we have
\begin{align*}
	&2\tau\| g^n\|_{L^2(\partial\Omega)}\|U^n(f)-U^{n-1}(f)\|_{L^2(\partial\Omega)} \\
	&~\quad \le C\tau  \| g\|_{L^\infty(t^{n-1},t^n; L^2(\partial\Omega))} \big(\|U^n(f)\|_{H^1(\Omega)} + \|U^{n-1}(f)\|_{H^1(\Omega)}\big)\\
	&~\quad \le \frac{C_1\tau}{2}\|U^n(f)\|^2_{H^1(\Omega)} + C\tau\left( \|U^{n-1}(f)\|^2_{H^1(\Omega)} + \|g\|^2_{L^\infty \left( t^{n-1}, t^n;L^2(\partial\Omega)\right) }\right)
	\end{align*}	
	We thus arrive at
	\begin{align}\label{6-2-17ct1}
	&\|U^n(f)-U^{n-1}(f)\|^2_{L^2(\Omega)} + \tau \|U^n(f)\|^2_{H^1(\Omega)} \notag\\
	&~\quad \le C\tau\left( \|U^{n-1}(f)\|^2_{H^1(\Omega)}+\int_{t^{n-1}}^{t^n}\|f\|^2_{L^2(\Omega)} dt+\|g\|^2_{L^\infty \left( t^{n-1}, t^n;L^2(\partial\Omega)\right) }\right)
	\end{align}
	as $\epsilon$ is small enough with $C$ independent of $n$. Since $\|U^0(f)\|_{H^1(\Omega)} = \|\Pi_hq\|_{H^1(\Omega)} \le C\|q\|_{H^1(\Omega)}$, we deduce from \cref{6-2-17ct1} that
	\begin{align*}
	&\|U^n(f)\|^2_{H^1(\Omega)} \\
	&~\quad \le C\left( \| q \|^2_{H^1(\Omega)}+ \|g\|^2_{L^\infty\left( 0,T; L^2(\partial\Omega)\right) }  + \|f\|^2_{L^2(\Omega_T)} \right) + \sum_{m=0}^{n-1} \alpha_m\|U^m(f)\|^2_{H^1(\Omega)}
	\end{align*}
for all $n\ge 1$, where $\alpha_0 =\ldots=\alpha_{n-2}=0$ and $\alpha_{n-1}=C$. Therefore, an application of the discrete Gronwall inequality implies \cref{19-1-17ct4}
	for all $n\ge 0$.
	
	(ii) By \cref{18-1-17ct5*}, for all $\theta\in L^2\left(0,T; H^1(\Omega) \right)$ we rewrite
	\begin{align}\label{6-2-17ct2}
	( \partial U^n(f),\theta)_{L^2(\Omega)} 
	&= (\partial U^n(f),\Pi_h \theta)_{L^2(\Omega)} + (\partial U^n(f),\theta-\Pi_h \theta)_{L^2(\Omega)} \notag\\
	&= (\overline{f}^n,\Pi_h \theta)_{L^2(\Omega)} +( g^n,\Pi_h \theta)_{L^2(\partial\Omega)} - a^n(\overline{\partial} U^n(f),\Pi_h \theta) \notag\\
	&~\quad + (\partial U^n(f),\theta-\Pi_h \theta)_{L^2(\Omega)}.
	\end{align}
	We have that
	\begin{align*}
	\sum_{n=1}^M \int_{t^{n-1}}^{t^n}(\overline{f}^n,\Pi_h \theta)_{L^2(\Omega)}dt
	&\le \sum_{n=1}^M \left( \int_{t^{n-1}}^{t^n}\|f\|^2_{L^2(\Omega)}dt\right)^{1/2}\left( \int_{t^{n-1}}^{t^n}\|\Pi_h\theta\|^2_{L^2(\Omega)}dt\right)^{1/2} \\
	&\le C \|f\|_{L^2(\Omega_T)}\|\theta\|_{L^2\left(0,T; H^1(\Omega)\right) }
	\end{align*}
	and
	\begin{align*}
	&\sum_{n=1}^M \int_{t^{n-1}}^{t^n}( g^n,\Pi_h \theta)_{L^2(\partial\Omega)} dt\\
	&~\quad \le \sum_{n=1}^M \left( \tau\|g\|^2_{L^\infty\left( t^{n-1}, t^n; L^2(\partial\Omega)\right) }\right)^{1/2}\left( \int_{t^{n-1}}^{t^n}\|\Pi_h\theta\|^2_{L^2(\partial\Omega)}dt\right)^{1/2} \\
	&~\quad \le C \|g\|_{L^\infty\left( 0,T; L^2(\partial\Omega)\right) }\|\theta\|_{L^2\left(0,T; H^1(\Omega)\right) }.
	\end{align*}
	Using \cref{19-1-17ct4}, we further get that
	\begin{align*}
	&\sum_{n=1}^M \int_{t^{n-1}}^{t^n}a^n(\overline{\partial} U^n(f),\Pi_h \theta)dt\\
	&~\quad \le C\sum_{n=1}^M \int_{t^{n-1}}^{t^n}\left( \|U^n(f)\|_{H^1(\Omega)}+\|U^{n-1}(f)\|_{H^1(\Omega)}\right) \|\Pi_h\theta\|_{H^1(\Omega)}dt \\
	&~\quad \le C\sum_{n=1}^M \int_{t^{n-1}}^{t^n}\|U^n(f)\|_{H^1(\Omega)}\|\theta\|_{H^1(\Omega)} dt\\
	&~\quad \le C\sum_{n=1}^M \left( \tau\|U^n(f)\|^2_{H^1(\Omega)}\right)^{1/2}\left( \int_{t^{n-1}}^{t^n}\|\theta\|^2_{H^1(\Omega)} dt\right)^{1/2}\\
	&~\quad \le C\left( \sum_{n=1}^M \tau\|U^n(f)\|^2_{H^1(\Omega)}\right)^{1/2}\left( \sum_{n=1}^M \int_{t^{n-1}}^{t^n}\|\theta\|^2_{H^1(\Omega)} dt\right)^{1/2}\\
	&~\quad \le C\left( \| q \|_{H^1(\Omega)}+ \|g\|_{L^\infty\left( 0,T; L^2(\partial\Omega)\right) }  + \|f\|_{L^2(\Omega_T)} \right)\|\theta\|_{L^2\left(0,T; H^1(\Omega)\right) }.
	\end{align*}
	Next, we deduce from \cref{6-2-17ct1} that
	\begin{align}\label{6-2-17ct3}
	&\sum_{n=1}^M\|U^n(f)-U^{n-1}(f)\|^2_{L^2(\Omega)} \notag\\
	&~\quad \le C\left( \| q \|^2_{H^1(\Omega)}+ \|g\|^2_{L^\infty\left( 0,T; L^2(\partial\Omega)\right) }  + \|f\|^2_{L^2(\Omega_T)} \right).
	\end{align}
	Utilizing \cref{23/5:ct1}, we have that
	\begin{align*}
	&\sum_{n=1}^M \int_{t^{n-1}}^{t^n}(\partial U^n(f),\theta-\Pi_h \theta)_{L^2(\Omega)} dt\\
	&~\quad \le  \left( \sum_{n=1}^M \tau\|\partial U^n(f)\|^2_{L^2(\Omega)}\right)^{1/2}\left( \sum_{n=1}^M\int_{t^{n-1}}^{t^n}\|\theta-\Pi_h\theta\|^2_{L^2(\Omega)} dt\right)^{1/2}\notag\\
	&~\quad\le C\left( \sum_{n=1}^M h^2\tau\|\partial U^n(f)\|^2_{L^2(\Omega)}\right)^{1/2}\|\theta\|_{L^2\left(0,T; H^1(\Omega)\right) }\notag\\
	&~\quad \le C\left( \| q \|_{H^1(\Omega)}+ \|g\|_{L^\infty\left( 0,T; L^2(\partial\Omega)\right) }  + \|f\|_{L^2(\Omega_T)} \right)\|\theta\|_{L^2\left(0,T; H^1(\Omega)\right) },
	\end{align*}
	by \cref{6-2-17ct3} as $h^2\tau^{-1} \le C$. Therefore, \cref{19-1-17ct4*} follows from \cref{6-2-17ct2} and the above estimates. Furthermore, in the same manner we also get \cref{19-1-17ct4***}.
	
	(iii) We consider the piecewise constant function with respect to $t$ defined as follows
	$${\Phi_{h,\tau}}_{|\Omega\times (t^{n-1},t^n]} := U^n(f)-U^{n-1}(f) \quad\mbox{for all}\quad n\in I.$$
	Due to \cref{19-1-17ct4}, the sequence $\left( \Phi_{h,\tau}\right)_{h,\tau} $ is bounded in the $L^2\left(0,T; H^1(\Omega)\right)$-norm. Therefore, there exist a subsequence of it denoted the same symbol and an element $\Phi\in L^2\left(0,T; H^1(\Omega)\right)$ such that for all $\theta\in L^2\left(0,T; H^1(\Omega) \right)$
	\begin{align*}
	\lim_{h,\tau \to 0}\int_0^T\left( \Phi_{h,\tau},\theta\right)_{H^1(\Omega)} dt = \int_0^T\left( \Phi,\theta\right)_{H^1(\Omega)} dt
	\end{align*}
and $\|\Phi\|_{L^2(\Omega_T)} \le \liminf_{h,\tau \to 0}\left\| \Phi_{h,\tau}\right\|_{L^2(\Omega_T)}$, since the mapping $L^2\left(0,T; H^1(\Omega)\right) \hookrightarrow L^2(\Omega_T)$ is continuous. On the other hand, it follows from \cref{6-2-17ct3} that
	\begin{align*}
	\left\| \Phi_{h,\tau}\right\|^2_{L^2(\Omega_T)} = \sum_{n=1}^M \int_{t^{n-1}}^{t^n}\left\| \Phi_{h,\tau}\right\|^2_{L^2(\Omega)} = \sum_{n=1}^M\tau\|U^n(f)-U^{n-1}(f)\|^2_{L^2(\Omega)} \to 0
	\end{align*}
	as $h,\tau\to 0$, which finishes the proof.
\end{proof}

\begin{lemma}\label{conv-dis}
	Assume that the sequence $\left( f_k\right)$ weakly
	converges  in $L^2(\Omega_T)$ to an element $f$. Then for any fixed $n\in I_0$, the sequence $\left(U^n(f_k)\right)$ converges to $U^n(f)$ in the $H^1(\Omega)$-norm.
\end{lemma}

\begin{proof}
	The proof is based on standard arguments, it is therefore omitted here.
\end{proof}

\begin{theorem}\label{disc.optim}
	The problem $\left(\mathcal{P}_{\rho,\delta,h,\tau}\right)$
	attains a unique solution $f:= f_{\rho,\delta,h,\tau}$ satisfying the equation
	\begin{align}\label{23-1-17ct2}
	f_{|\Omega\times(t^{n-1},t^n]} = f^*- \rho^{-1}P^{n-1}(f)
	\end{align}
	for any $n\in I$.
\end{theorem}

\begin{proof}
	In virtue of \Cref{conv-dis}, it is straightforward to show the uniqueness existence of a solution $f$ to $\left(\mathcal{P}_{\rho,\delta,h,\tau}\right)$. We now show the equation \cref{23-1-17ct2}. We have $J'_{\rho,\delta,h,\tau}(f)\xi =0$ for all $\xi\in L^2(\Omega_T)$.
	By \cref{18-1-17ct5**}--\cref{18-1-17ct5***}, we get
	\begin{align*}
	J'_{\rho,\delta,h,\tau}(f)\xi
	&= 2\sum_{n=1}^M \int_{t^{n-1}}^{t^n}\left( U^n(f)-z_\delta, {U^n}'(f)\xi\right) _{L^2(\Gamma)}dt + 2\rho(f-f^*,\xi)_{L^2(\Omega_T)}\\
	&= 2\sum_{n=1}^M \int_{t^{n-1}}^{t^n}\left( U^n(f)-z_\delta, \widehat{U}^n(\xi)\right)_{L^2(\Gamma)}dt + 2\rho(f-f^*,\xi)_{L^2(\Omega_T)}\\
	&= 2\tau\sum_{n=1}^M-\left( \partial P^n(f), \widehat{U}^n(\xi) \right)_{L^2(\Omega)} + 2\tau\sum_{n=1}^Ma^n\left( \overline{\partial}P^n(f),\widehat{U}^n(\xi)\right) \\
	&~\quad + 2\rho(f-f^*,\xi)_{L^2(\Omega_T)}.
	\end{align*}
	Using the identities
	\begin{equation}\label{26-1-17ct1}
	\begin{aligned}
	&\sum_{n=1}^M (\alpha^n-\alpha^{n-1})\beta^n = \alpha^M\beta^M-\alpha^0\beta^0 -\sum_{n=1}^M (\beta^n-\beta^{n-1})\alpha^{n-1},\\
	&\sum_{n=1}^M (\alpha^n+\alpha^{n-1})\beta^n = \alpha^M\beta^M-\alpha^0\beta^0 +\sum_{n=1}^M (\beta^n+\beta^{n-1})\alpha^{n-1}
	\end{aligned}
	\end{equation}
	together with $P^M(f)=\widehat{U}^0(\xi)=0$, we obtain
	\begin{align}\label{25-3-19ct1}
	J'_{\rho,\delta,h,\tau}(f)\xi
	&= 2\tau\sum_{n=1}^M\left( \partial\widehat{U}^n(\xi), P^{n-1}(f) \right)_{L^2(\Omega)} + 2\tau\sum_{n=1}^Ma^n\left( \overline{\partial}\widehat{U}^n(\xi),P^{n-1}(f)\right) \notag \\
	&~\quad + 2\rho(f-f^*,\xi)_{L^2(\Omega_T)}\notag\\
	&=2\tau\sum_{n=1}^M(P^{n-1}(f),\overline{\xi}^n)_{L^2(\Omega)}+ 2\rho(f-f^*,\xi)_{L^2(\Omega_T)}.
	\end{align}
	For any fixed $n\in I$ and for all $\varphi\in L^2(\Omega\times(t^{n-1},t^n])$ we consider $\xi := \varphi\chi_{\Omega \times (t^{n-1},t^n]}\in L^2(\Omega_T)$ and then have $\overline{\xi}^k = 0$ as $k\neq n$ and
	\begin{align*}
	\tau\sum_{k=1}^M(P^{k-1}(f),\overline{\xi}^k)_{L^2(\Omega)}
	&= \tau(P^{n-1}(f),\overline{\xi}^n)_{L^2(\Omega)} \\
	&=\tau \int_\Omega \left( P^{n-1}(x;f)\frac{1}{\tau}\int_{t^{n-1}}^{t^n}\varphi(x,t)dt\right)dx \\
	&= \int_{t^{n-1}}^{t^n}(P^{n-1}(f),\varphi)_{L^2(\Omega)}dt
	\end{align*}
	as well as
	$$(f-f^*,\xi)_{L^2(\Omega_T)} =\int_{t^{n-1}}^{t^n}(f-f^*,\varphi)_{L^2(\Omega)}dt.$$
	Thus, we arrive at
	\begin{align*}
	\int_{t^{n-1}}^{t^n}\left( P^{n-1}(f) + \rho(f-f^*),\varphi\right) _{L^2(\Omega)}dt =0,
	\end{align*}
	where $\varphi\in L^2(\Omega\times(t^{n-1},t^n])$ is arbitrary. This implies \cref{23-1-17ct2}. The proof is finished.
\end{proof}

\begin{remark}
	For any fixed $f\in L^2(\Omega_T)$, denote by
	$${\mathbb{G}_J(x,t;f)}_{|\Omega\times [t^{n-1},t^n)} := P^{n-1}(x,t;f) \quad\mbox{with}\quad n\in I.$$
	In view of the identity \cref{25-3-19ct1}, the $L^2$-gradient of the cost functional at $f$ is given by
	\begin{align}\label{25-3-19ct2}
	\nabla J_{\rho,\delta,h,\tau}(x,t;f) = 2\mathbb{G}_J(x,t;f) +2 \rho(f-f^*)
	\end{align}
	i.e. the equality
	$$J'_{\rho,\delta,h,\tau}(f)\xi = (\nabla J_{\rho,\delta,h,\tau}(f),\xi)_{L^2(\Omega_T)}$$
	holds true for all $\xi\in L^2(\Omega_T)$.
\end{remark}

\section{Convergence of finite dimensional approximations} \label{convergence}

The aim of this section is to show finite dimensional approximations, i.e. solutions of $\left(\mathcal{P}_{\rho,\delta,h,\tau}\right)$, converge to the sought source. To do so, we state some auxiliary results.

\begin{lemma}\label{limit}
	(i) For all $\phi\in L^2\left( 0,T;H^k(\Omega)\right)$ with $k \in \{0,1\}$ there hold
	\begin{equation}\label{23/10:ct2*}
	\begin{aligned}
	&\lim_{h\to 0} \left\| \phi - \Pi_h \phi
	\right\|_{L^2\left( 0,T;H^k(\Omega)\right) } =0,\quad \lim_{\tau\to 0} \left\| \phi -  \phi_\tau
	\right\|_{L^2\left( 0,T;H^k(\Omega)\right) } =0 \quad\mbox{and}\\
	&\lim_{h,\tau\to 0} \big\| \phi -\Pi_h  \phi_\tau
	\big\|_{L^2\left( 0,T;H^k(\Omega)\right) } =0.
	\end{aligned}
	\end{equation}
	(ii) Assume that $\phi\in L^\infty\left( \Omega_T\right)$, then
	\begin{align}\label{29-1-17ct1}
	\lim_{\tau\to 0} \int_0^T\int_\Omega ( \phi_\tau-\phi) uv_\tau dxdt=0
	\end{align}
	for all $u\in L^2(\Omega_T)$ and any bounded sequence $(v_\tau)_\tau$ in $L^2(\Omega_T)$.
\end{lemma}

\begin{proof}
	(i) The first and second statements of \cref{23/10:ct2*} follow directly from \cref{23/10:ct2} and \cref{1-2-17ct1*}, resectively, and the Lebesgue's dominated convergence theorem. Meanwhile, for a.e in $t\in (0,T)$, since
	\begin{align*}
	\big\| \phi-\Pi_h \phi_\tau\big\|_{H^k(\Omega)} 
	&\le \big\| \phi-\Pi_h\phi\big\|_{H^k(\Omega)} + \big\| \Pi_h(\phi- \phi_\tau)\big\|_{H^k(\Omega)} \\
	&\le \big\| \phi-\Pi_h\phi\big\|_{H^k(\Omega)} + C\big\| \phi- \phi_\tau\big\|_{H^k(\Omega)},
	\end{align*}
	the third assertion follows from the first and second ones.
	
	For (ii) we take an arbitrary $\epsilon>0$ and $u^\epsilon\in C(\overline{\Omega_T})$ such that $\|u-u^\epsilon\|_{L^2(\Omega_T)} <\epsilon$. Then we have
	\begin{align*}
	&\left|\int_0^T\int_\Omega ( \phi_\tau-\phi) uv_\tau dxdt\right| \\
	&~\quad \le \int_0^T\int_\Omega | \phi_\tau-\phi| |u^\epsilon||v_\tau| dxdt + \int_0^T\int_\Omega | \phi_\tau-\phi| |u-u^\epsilon||v_\tau| dxdt\\
	&~\quad \le \|u^\epsilon\|_{C(\overline{\Omega_T})} \|v_\tau\|_{L^2(\Omega_T)}\| \phi_\tau-\phi \|_{L^2(\Omega_T)} + \left( 2\|\phi\|_{L^\infty\left( \Omega_T\right)}\|v_\tau\|_{L^2(\Omega_T)}\right)\epsilon.
	\end{align*}
	Sending $\tau$ to zero, we thus have that
	$
	\lim_{\tau\to 0}\left|\int_0^T\int_\Omega ( \phi_\tau-\phi) uv_\tau dxdt\right| \le C\epsilon
	$
	for all $\epsilon>0$ with the constant $C$ independent of $\epsilon$, which yields \cref{29-1-17ct1}. The proof is completed.
\end{proof}

\begin{lemma}\label{6-5-16ct6}
	Assume that $\limk h_k = \limk \tau_k = 0$ and the sequence $\left( z_{\delta_k}\right) \subset L^2(\Sigma)$ converge to $z_\delta $ in the  $L^2(\Sigma)$-norm.
	Let the sequence $(f_k) \subset L^2(\Omega_T)$  weakly converge  in $L^2(\Omega_T)$ to $f$. Then,
	\begin{align} \label{20-5-16ct3}
	 \liminf_{k\to\infty} \sum_{n=1}^{M_k} \int_{t^{n-1}}^{t^n}\|U^n_{h_k,\tau_k}(f_k)-z_{\delta_k}\|^2_{L^2(\Gamma)}dt \ge \|u(f)-z_\delta\|^2_{L^2(\Sigma)},
	\end{align}
	where $M_k =T/\tau_k$ and $U^n_{h_k,\tau_k}(f_k)$ defined by \cref{18-1-17ct5*}.
\end{lemma}

\begin{proof}
	For convenience of exposition we denote by $U^n_k := U^n_{h_k,\tau_k}(f_k)$. Let $\Phi_k =\Phi_k(x,t)$ be the piecewise linear, continuous interpolation of $\left( U^n_k\right)_{n= 0,\ldots,M_k}$ with respect to $t$, i.e.
	$$\Phi_k(x,t) := (t-t^{n-1})\partial U^n_k + U^{n-1}_k$$
	with $\partial U^n_k = \tau_k^{-1}(U^n_k-U^{n-1}_k)$ and $(x,t)\in\Omega\times(t^{n-1},t^n]$, $n = 1,\ldots, M_k$. We first note that for all $t\in (t^{n-1},t^n)$
	\begin{align}\label{1-4-19ct1}
	\frac{\partial \Phi_k}{\partial t} = \partial U^n_k \quad \mbox{and} \quad
	\int_{t^{n-1}}^{t^n}\Phi_k dt = \tau_k \overline{\partial}U^n_k.
	\end{align}
	Thus, $\left( \Phi_k\right)\subset H^1\left(0,T; \mathcal{V}^1_h\right) \subset H^1\left(0,T; H^1(\Omega)\right) \subset \mathcal{W}(0,T)$. Further, the inequalities \cref{19-1-17ct4} and \cref{19-1-17ct4***} yield that the sequence $\left( \Phi_k\right)$ is bounded in the reflexive space $\mathcal{W}(0,T)$. There exists a subsequence of $\left( \Phi_k\right)$ denoted again by $\left( \Phi_k\right)$ and an element $u\in \mathcal{W}(0,T)$ such that $\left( \Phi_k\right)$ weakly converges in $\mathcal{W}(0,T)$ to $u$.
	
	We show that $u=u(f)$. In fact, for all $\phi\in L^2\left( 0,T, H^1(\Omega)\right) $ we have
	\begin{align}\label{31-1-17ct2}
	&\int_0^T\left\langle \frac{\partial \Phi_k}{\partial t}, \phi \right\rangle_{\left( H^1(\Omega)^*, H^1(\Omega)\right) }dt + \int_0^T a(\Phi_k,\phi)dt\notag\\
	&~\quad = \sum_{n=1}^{M_k} \int_{t^{n-1}}^{t^n} \left( \partial U^n_k,\phi\right)_{L^2(\Omega)}dt + \int_0^T\int_\Omega  A\nabla \Phi_k\cdot\nabla \phi dxdt \notag\\
	&~\quad + \int_0^T\int_\Omega
	b\Phi_k\phi dxdt + \int_0^T\int_{\partial\Omega}\sigma\Phi_k\phi dxdt.
	\end{align}
	For all $t\in (t^{n-1},t^n]$, we have
	\begin{align}\label{31-1-17ct2*}
	\big( \partial U^n_k,\phi\big)_{L^2(\Omega)} = \big( \partial U^n_k,\Pi_{h_k} \overline{\phi}^n\big)_{L^2(\Omega)} + \big( \partial U^n_k, \phi - \Pi_{h_k} \overline{\phi}^n\big)_{L^2(\Omega)}
	\end{align}
	and, by \cref{19-1-17ct4*} and \cref{23/10:ct2*},
	\begin{align}\label{31-1-17ct2**}
	\left| \sum_{n=1}^{M_k} \int_{t^{n-1}}^{t^n} \big( \partial U^n_k, \phi - \Pi_{h_k} \overline{\phi}^n\big)_{L^2(\Omega)}\right| \le C\big\| \phi - \Pi_{h_k} \phi_{\tau_k}\big\|_{L^2\left( 0,T; H^1(\Omega)\right) } \to 0
	\end{align}
	as $k\to\infty$. Using the estimates  \cref{19-1-17ct4} and \cref{19-1-17ct4*} as well as the equalities \cref{23/10:ct2*}, \cref{29-1-17ct1},
	we decompose the remainder in the right hand side of \cref{31-1-17ct2} as follows
	\begin{align}\label{7-2-17ct4}
	&\int_0^T\int_\Omega  A\nabla \Phi_k\cdot\nabla \phi dxdt + \int_0^T\int_\Omega b\Phi_k\phi dxdt + \int_0^T\int_{\partial\Omega}\sigma\Phi_k\phi dxdt \notag\\
	&~\quad  = \sum_{n=1}^{M_k} \int_{t^{n-1}}^{t^n}\int_\Omega  A^n\nabla \Phi_k\cdot\nabla \Pi_{h_k} \overline{\phi}^n dxdt + \sum_{n=1}^{M_k} \int_{t^{n-1}}^{t^n}\int_\Omega
	b^n\Phi_k\Pi_{h_k} \overline{\phi}^n dxdt \notag\\
	&~\quad\qquad + \sum_{n=1}^{M_k} \int_{t^{n-1}}^{t^n}\int_{\partial\Omega}\sigma^n \Phi_k\Pi_{h_k} \overline{\phi}^n dxdt + \mathcal{R}_k\notag\\
	&~\quad = \sum_{n=1}^{M_k} \tau_k a^n \big(\overline{\partial}U^n_k,\Pi_{h_k} \overline{\phi}^n\big)+ \mathcal{R}_k,
	\end{align}
	by \cref{1-4-19ct1}, where
	\begin{align*}
	&\mathcal{R}_k\\
	&~:= \sum_{n=1}^{M_k} \int_{t^{n-1}}^{t^n}\int_\Omega (A-A^n) \nabla \Phi_k\cdot\nabla \phi dxdt + \sum_{n=1}^{M_k} \int_{t^{n-1}}^{t^n}\int_\Omega  A^n\nabla \Phi_k\cdot\nabla  (\phi - \Pi_{h_k}\overline{\phi}^n) dxdt\\
	&~+ \sum_{n=1}^{M_k} \int_{t^{n-1}}^{t^n}\int_\Omega (b-b^n) \nabla \Phi_k\cdot\nabla \phi dxdt + \sum_{n=1}^{M_k} \int_{t^{n-1}}^{t^n}\int_\Omega  b^n\nabla \Phi_k\cdot\nabla  (\phi - \Pi_{h_k}\overline{\phi}^n) dxdt\\
	&~+ \sum_{n=1}^{M_k} \int_{t^{n-1}}^{t^n}\int_{\partial\Omega} (\sigma-\sigma^n) \nabla \Phi_k\cdot\nabla \phi dxdt + \sum_{n=1}^{M_k} \int_{t^{n-1}}^{t^n}\int_{\partial\Omega}  \sigma^n\nabla \Phi_k\cdot\nabla  (\phi - \Pi_{h_k}\overline{\phi}^n) dxdt.
	\end{align*}
	We remark that due to the continuity of data and \cref{23/10:ct2*}, the relation
	$
	\limk \mathcal{R}_k=0
	$
	holds true. Therefore, we obtain from \cref{31-1-17ct2}--\cref{7-2-17ct4} that
	\begin{align*}
	&\int_0^T \left\langle \frac{\partial u}{\partial t}, \phi \right\rangle_{\left( H^1(\Omega)^*, H^1(\Omega)\right) } dt + \int_0^T a(u,\phi) dt\\
	&~\quad = \limk\left( \int_0^T\left\langle \frac{\partial \Phi_k}{\partial t}, \phi \right\rangle_{\left( H^1(\Omega)^*, H^1(\Omega)\right) }dt + \int_0^T a(\Phi_k,\phi)dt\right) \notag\\
	&~\quad = \limk \sum_{n=1}^{M_k} \tau_k\left(  \big( \partial U^n_k,\Pi_{h_k} \overline{\phi}^n\big)_{L^2(\Omega)} + a^n \big(\overline{\partial}U^n_k,\Pi_{h_k} \overline{\phi}^n\big) \right) \notag\\
	&~\quad = \limk \sum_{n=1}^{M_k} \tau_k \left( \big( \overline{f_k}^n,\Pi_{h_k} \overline{\phi}^n\big)_{L^2(\Omega)} + \big(   g^n, \Pi_{h_k} \overline{\phi}^n\big)_{L^2(\partial\Omega)}\right)\\
	&~\quad = (f,\phi)_{L^2(0,T; L^2(\Omega))} + ( g, \phi)_{L^2(0,T;L^2(\partial\Omega))},
	\end{align*}
	here we used \cref{18-1-17ct5*}. Further, with standard arguments it can be shown that $u(x,0) = q(x)$. Thus, we get that the sequence $\left( \Phi_k\right)$ weakly converges  in $\mathcal{W}(0,T)$ to $u(f)$.
	
	Next, for each $\theta\in C^1\left( [0,T]; H^1(\Omega)\right) $ we will show that
	\begin{align}\label{7-2-17ct5}
	\limk \sum_{n=1}^{M_k} \int_{t^{n-1}}^{t^n} \big(\Phi_k-U^n_k, \theta\big)_{H^1(\Omega)}dt =0,
	\end{align}
	which then holds true for each $\theta\in L^2\left( 0,T; H^1(\Omega)\right)$, by the density argumentation. We have that
	\begin{align*}
	&\int_{t^{n-1}}^{t^n} \big(\Phi_k-U^n_k, \theta\big)_{H^1(\Omega)}dt \\
	&~\quad = \frac{1}{2}\int_{t^{n-1}}^{t^n} \big(\partial U^n_k, \theta\big)_{H^1(\Omega)}d(t-t^{n-1}-\tau_k)^2 \\
	&~\quad = \frac{\tau_k^2}{2}\big(\partial U^n_k, \theta(t^{n-1})\big)_{H^1(\Omega)} - \frac{1}{2}\int_{t^{n-1}}^{t^n} (t-t^{n-1}-\tau_k)^2\left(\partial U^n_k, \frac{d\theta}{dt}\right)_{H^1(\Omega)}dt
	\end{align*}
	and so
	\begin{align*}
	&\left| \sum_{n=1}^{M_k}\int_{t^{n-1}}^{t^n} \big(\Phi_k-U^n_k, \theta\big)_{H^1(\Omega)}dt \right|\\
	&~\quad \le \frac{1}{2}\sum_{n=1}^{M_k}\int_{t^{n-1}}^{t^n}\big(U^n_k-U^{n-1}_k, \theta(t^{n-1})\big)_{H^1(\Omega)} dt \\
	&~\qquad +C\tau_k\|\theta\|_{C^1\left( [0,T]; H^1(\Omega)\right) }\sum_{n=1}^{M_k}\tau^2_k \|\partial U^n_k\|_{H^1(\Omega)} \\
	&~\quad \to 0 \quad\mbox{as}\quad k\to\infty,
	\end{align*}
	by \cref{19-1-17ct4**} and \cref{19-1-17ct4}.
	
	Finally, using the continuity of the mapping $ L^2\left(0,T; H^1(\Omega) \right) \to L^2\left( 0,T;L^2(\Gamma)\right)$ and the fact $\|z_{\delta_k} - z_\delta\|_{L^2(\Sigma)} \to 0$ as $k\to \infty$, we arrive at \cref{20-5-16ct3}, which finishes the proof.
\end{proof}

Next we introduce the notion of the $f^*$-minimum-norm solution of the identification problem.

\begin{lemma}\label{identification_problem}
	The problem
	$$
	\min_{f \in \mathcal{I}\left(Z\right)}  \|f-f^*\|_{L^2(\Omega_T)} \eqno\left(\mathcal{IP}\right)
	$$
	attains a unique solution, which is called the {\it $f^*$-minimum-norm solution} of the identification problem, where
	\begin{align}\label{3-5-16ct1}
	\mathcal{I} \left(Z\right) :=\left\{ f \in L^2(\Omega_T) ~\big|~ u(f)_{|\Sigma} =  Z \right\}.
	\end{align}
\end{lemma}

\begin{proof}
	Due to \Cref{weakly conv.}, $\mathcal{I} \left(Z\right) \neq \emptyset$ is a close subset of $L^2(\Omega_T)$. Furthermore, it is a convex set.
	Therefore, the minimization problem has a unique solution, which finishes the proof.
\end{proof}

We now show the main result of this section on the convergence of finite dimensional approximations $f_{\rho,\delta,h,\tau}$ of $\left(\mathcal{P}_{\rho,\delta,h,\tau}\right)$ to the $f^*$-minimizing-norm solution of the idendification problem $\left(\mathcal{IP}\right)$.
For any fixed $f\in L^2(\Omega_T)$, let $u(f)$ and $U^n(f)$ define by \cref{17-10-16ct2} and \cref{18-1-17ct5*}, respectively. We recall the convergence of the Crank-Nicolson Galerkin method for linear parabolic problems
\begin{align}\label{4-3-19ct1}
\lim_{h,\tau\to 0} \omega_{h,\tau}(f) = 0 \quad\mbox{with}\quad \omega_{h,\tau}(f) := \sum_{n=1}^M \int_{t^{n-1}}^{t^n}\left\| U^n(f) -u(f)\right\|^2_{H^1(\Omega)} dt.
\end{align}

\begin{theorem}\label{stability2}
	Let $f^\dag$ be the unique $f^*$-minimum-norm solution of the identification problem $\left(\mathcal{IP}\right)$. Let $(h_k), (\tau_k)$, $(\delta_k)$ and $(\rho_k)$ be any positive sequences such that
	\begin{align}\label{31-8-16ct1}
	h_k\to 0, \quad \tau_k \to 0, \quad \rho_k \to 0, \quad \frac{\delta_k^2}{\rho_k} \to 0, \quad \frac{\omega_{h_k,\tau_k}(f^\dag)}{\rho_k} \to 0 
	\end{align}
	as $k\to\infty$. Furthermore, assume that
	$\left( z_{\delta_k}\right) \subset L^2(\Sigma)$ is a sequence satisfying
	$$\|z_{\delta_k} - Z\|_{L^2(\Sigma)} \le \delta_k$$
	and $f_k$ denotes the unique minimizer of $(\mathcal{P}_{\rho_k,\delta_k, h_k, \tau_k})$ for each $k\in \mathbb{N}$. Then:
	
	(i) The sequence $(f_k)$ converges to $f^\dag$ in the $L^2(\Omega_T)$-norm.
	
	(ii) The following equality holds true
	$$\limk\sum_{n=1}^{M_k}\int_{t^{n-1}}^{t^n}\|U^n_{h_k,\tau_k}(f_k)-u(f^\dag)\|^2_{H^1(\Omega)}dt =0.$$
\end{theorem}

\begin{proof}
	For $n= 1,\ldots,M_k$ we write $U^n_k$ and ${U^n_k}^\dag$ instead $U^n_{h_k,\tau_k}(f_k)$ and $U^n_{h_k,\tau_k}(f^\dag)$ for short, respectively. Since $f_k$ is the solution of $(\mathcal{P}_{\rho_k,\delta_k, h_k, \tau_k})$, we have 
	\begin{align}\label{6-5-16ct8}
	&\sum_{n=1}^{M_k} \int_{t^{n-1}}^{t^n}\|U^n_k-z_{\delta_k}\|^2_{L^2(\Gamma)} dt + \rho_k\|f_k-f^*\|^2_{L^2(\Omega_T)} \notag\\
	&~\quad \le \sum_{n=1}^{M_k} \int_{t^{n-1}}^{t^n}\|{U^n_k}^\dag-z_{\delta_k}\|^2_{L^2(\Gamma)}dt + \rho_k\|f^\dag-f^*\|^2_{L^2(\Omega_T)}.
	\end{align}
	We get
	\begin{align}\label{1-2-17ct4}
	&\sum_{n=1}^{M_k} \int_{t^{n-1}}^{t^n}\|{U^n_k}^\dag-z_{\delta_k}\|^2_{L^2(\Gamma)} dt
	= \sum_{n=1}^{M_k} \int_{t^{n-1}}^{t^n}\|{U^n_k}^\dag- u(f^\dag) + Z -z_{\delta_k}\|^2_{L^2(\Gamma)} dt \notag\\
	&~\quad \le 2\sum_{n=1}^{M_k} \int_{t^{n-1}}^{t^n}\|{U^n_k}^\dag-u(f^\dag)\|^2_{L^2(\Gamma)}dt +2\sum_{n=1}^{M_k} \int_{t^{n-1}}^{t^n}\|Z-z_{\delta_k}\|^2_{L^2(\Gamma)}dt \notag\\
	&~\quad \le 2\left( \omega_{h_k,\tau_k}(f^\dag) + \delta_k^2\right).
	\end{align}
	It follows from \cref{6-5-16ct8} and \cref{1-2-17ct4} that
	\begin{align}\label{6-5-16ct8*}
	&\sum_{n=1}^{M_k} \int_{t^{n-1}}^{t^n}\|U^n_k-z_{\delta_k}\|^2_{L^2(\Gamma)}dt + \rho_k\|f_k-f^*\|^2_{L^2(\Omega_T)}\notag\\
	&~\quad \le 2\left( \omega_{h_k,\tau_k}(f^\dag) +\delta_k^2\right)   + \rho_k\|f^\dag-f^*\|^2_{L^2(\Omega_T)}.
	\end{align}
	Therefore, by \cref{31-8-16ct1}, we have
	\begin{align}\label{6-5-16ct10}
	\limk \sum_{n=1}^{M_k} \int_{t^{n-1}}^{t^n} \|U^n_k-z_{\delta_k}\|^2_{L^2(\Gamma)} dt =0
	\end{align}
	and
	\begin{align}\label{6-5-16ct11}
	\limsupk \|f_k-f^*\|_{L^2(\Omega_T)} \le  \|f^\dag-f^*\|_{L^2(\Omega_T)}.
	\end{align}
	Applying \Cref{6-5-16ct6}, we deduce from the boundedness of $(f_k)$ due to \cref{6-5-16ct11} that there are a subsequence of it denoted by the same symbol and an element $\widehat{f} \in L^2(\Omega)$ such that
	\begin{equation}\label{6-5-16ct12}
	\begin{aligned}
	& f_k -f^* \rightharpoonup \widehat{f}-f^* \quad \mbox{weakly in}\quad L^2(\Omega_T) \\
	& \liminfk \|f_k-f^*\|_{L^2(\Omega_T)} \ge \|\widehat{f}-f^*\|_{L^2(\Omega_T)} \\
	&\liminf_{k\to\infty} \sum_{n=1}^{M_k} \int_{t^{n-1}}^{t^n}\|U^n_k - z_{\delta_k}\|^2_{L^2(\Gamma)}dt \ge \|u(\widehat{f})-Z\|^2_{L^2(\Sigma)}.
	\end{aligned}
	\end{equation}
	We thus obtain from \cref{6-5-16ct10} and \cref{6-5-16ct12} that
	\begin{align*}
	u(\widehat{f})_{|\Sigma}=Z \quad\mbox{or}\quad \widehat{f}\in \mathcal{I}(Z).
	\end{align*}
	Further, combining \cref{6-5-16ct11} with \cref{6-5-16ct12}, we also obtain
	\begin{align*}
	\|\widehat{f}-f^*\|_{L^2(\Omega_T)} 
	&\le \liminfk \|f_k-f^*\|_{L^2(\Omega_T)} \\
	&\le \limsupk \|f_k-f^*\|_{L^2(\Omega_T)} \le  \|f^\dag-f^*\|_{L^2(\Omega_T)}
	\end{align*}
	and so, by the uniqueness of the $f^*$-minimum-norm solution $f^\dag$,
	$$\widehat{f}=f^\dag \quad\mbox{and}\quad \limk \|f_k-f^\dag\|_{L^2(\Omega_T)} =0.$$
	Next, by \cref{18-1-17ct5*}, for all $\varphi_{h_k}\in \mathcal{V}^1_{h_k}$ and $n\in I$
	we have
	\begin{align}\label{4-3-19ct2}
	\big( \partial U^n_k-\partial {U^n_k}^\dag, \varphi_{h_k} \big)_{L^2(\Omega)} + a^n\big(\overline{\partial}U^n_k - \overline{\partial}{U^n_k}^\dag,\varphi_{h_k}\big) = \big(\overline{f_k}^n -\overline{f^\dag}^n,\varphi_{h_k}\big)_{L^2(\Omega)}.
	\end{align}
	Denoting by $e^n_k := U^n_k-{U^n_k}^\dag$ and taking $\varphi_{h_k} = 2\tau^2_k\partial e^n_k$ in the above equation, we can estimate the left hand side of \cref{4-3-19ct2} by
	\begin{align}\label{4-3-19ct3}
	&\big( \partial U^n_k-\partial {U^n_k}^\dag, \varphi_{h_k} \big)_{L^2(\Omega)} + a^n\big(\overline{\partial}U^n_k - \overline{\partial}{U^n_k}^\dag,\varphi_{h_k}\big) \notag\\
	&~\quad = 2\left(e^n_k - e^{n-1}_k, e^n_k - e^{n-1}_k\right)_{L^2(\Omega)} + \tau_k a^n\left( e^n_k + e^{n-1}_k, e^n_k - e^{n-1}_k\right)\notag\\
	&~\quad \ge 2\left\|e^n_k - e^{n-1}_k\right\|^2_{L^2(\Omega)} + C_1\tau_k\left\|e^n_k\right\|^2_{H^1(\Omega)} - C_2\tau_k\left\|e^{n-1}_k\right\|^2_{H^1(\Omega)}
	\end{align}
	and the right hand side by
	\begin{align}\label{4-3-19ct4}
	&\big(\overline{f_k}^n -\overline{f^\dag}^n,\varphi_{h_k}\big)_{L^2(\Omega)}
	= \int_\Omega \left( \frac{1}{\tau_k} \int_{t^{n-1}}^{t^n} (f_k(t) - f^\dag(t)) dt\right) \left( 2\tau_k (e^n_k - e^{n-1}_k)\right)dx \notag\\
	&~\quad = 2\int_\Omega \left(\int_{t^{n-1}}^{t^n} (f_k(t) - f^\dag(t)) dt\right) \left(e^n_k - e^{n-1}_k\right)dx \notag\\
	&~\quad \le 2 \left( \int_\Omega \left(\int_{t^{n-1}}^{t^n} (f_k(t) - f^\dag(t)) dt\right)^2 dx\right)^{1/2} \left( \int_\Omega \left(e^n_k - e^{n-1}_k\right)^2 dx\right)^{1/2} \notag\\
	&~\quad \le \epsilon \left\|e^n_k - e^{n-1}_k\right\|^2_{L^2(\Omega)} + \frac{C\tau_k}{\epsilon} \int_\Omega \int_{t^{n-1}}^{t^n} (f_k(t) - f^\dag(t))^2 dt dx
	\end{align}
	for any $\epsilon>0$. We then get from \cref{4-3-19ct2}--\cref{4-3-19ct4}
	\begin{align*}
	&2\left\|e^n_k - e^{n-1}_k\right\|^2_{L^2(\Omega)} + C_1\tau_k\left\|e^n_k\right\|^2_{H^1(\Omega)} \\
	&~\quad \le \epsilon \left\|e^n_k - e^{n-1}_k\right\|^2_{L^2(\Omega)} + \frac{C\tau_k}{\epsilon} \int_\Omega \int_{t^{n-1}}^{t^n} (f_k(t) - f^\dag(t))^2 dt dx + C_2\tau_k\left\|e^{n-1}_k\right\|^2_{H^1(\Omega)}
	\end{align*}
	that yields
	\begin{align*}
	\|e^n_k\|^2_{H^1(\Omega)} \le \int_\Omega \int_{t^{n-1}}^{t^n} C(f_k(t) - f^\dag(t))^2 dt dx + \sum_{m=0}^{n-1} \alpha_m\|e^{m}_k\|^2_{H^1(\Omega)},
	\end{align*}
	where $\alpha_0 = \ldots = \alpha_{n-2} = 0$ and $\alpha_{n-1}=C$.
	Since $e^0_k = 0$, applying Gronwall's inequality, we obtain
	\begin{align*}
	\|e^n_k\|^2_{H^1(\Omega)} \le C\int_\Omega \int_{t^{n-1}}^{t^n}(f_k(t) - f^\dag(t))^2 dt dx + C\int_\Omega \int_{t^{n-2}}^{t^{n-1}}(f_k(t) - f^\dag(t))^2 dt dx
	\end{align*}
	and so that
	\begin{align*}
	\sum_{n=1}^{M_k}\tau_k\|e^n_k\|^2_{H^1(\Omega)}\le C\|f_k-f^\dag\|^2_{L^2(\Omega_T)} \to 0
	\end{align*}
	as $k\to\infty$. Therefore, we in view of \cref{4-3-19ct1} conclude that
	\begin{align*}
	&\sum_{n=1}^{M_k}\int_{t^{n-1}}^{t^n}\|U^n_k-u(f^\dag)\|^2_{H^1(\Omega)}dt \\
	&~\quad \le 2\sum_{n=1}^{M_k}\int_{t^{n-1}}^{t^n}\|U^n_k-{U^n_k}^\dag\|^2_{H^1(\Omega)}dt + 2\sum_{n=1}^{M_k}\int_{t^{n-1}}^{t^n}\|{U^n_k}^\dag-u(f^\dag)\|^2_{H^1(\Omega)}dt \to 0
	\end{align*}
	as $k\to\infty$, which finishes the proof.
\end{proof}

\section{Convergence rates}\label{rate}
To obtain convergence rates for the Tikhonov regularization, some additional assumptions are required. We here assume that the data appearing in the system \cref{ct1} are regular enough such that the following error bound of the Crank-Nicolson Galerkin method for linear parabolic problems is fulfilled, see, e.g., \cite{thome}.
\begin{lemma}\label{29-3-19ct2}
	Let $u(f^\dag)$ and $U^n(f^\dag)$ be the solutions of \cref{17-10-16ct2} and \cref{18-1-17ct5*}, respectively. Then, the estimate
	\begin{align}\label{29-3-19ct3}
	\|U^n(f^\dag) - u^n(f^\dag)\|_{H^s(\Omega)}
	&\le Ch^{2-s}\left( \|q\|_{H^2(\Omega)} + \int_0^{t^n} \|u_t(f^\dag)\|_{H^2(\Omega)}dt\right) \notag\\
	&~\quad + C\tau^2h^{-s}\int_0^{t^n} \left(\|u_{ttt}(f^\dag)\|_{L^2(\Omega)} + \|\Delta u_{tt}(f^\dag)\|_{L^2(\Omega)}\right) dt
	\end{align}
	holds true for all $n\in I$, where $s\in \{0,1\}$.
\end{lemma}

We note that, due to the regularity theory for the parabolic equations (see, e.g.,  \cite[Theorem 6, pp.\ 365]{evan}), there exist the high order derivatives $\Delta u_{tt}(f^\dag)$ and $u_{ttt}(f^\dag)$ if $f^\dag_{tt} \in L^2(0,T; L^2(\Omega))$ and $f^\dag_{ttt} \in L^2(0,T; H^{-2}(\Omega))$, respectively.

To obtain the convergence rates of the regularized solutions to the identification, some smooth assumptions on the sought source $f^\dag$ and the exact data $u(f^\dag)$ should be addressed. We first assume they are such that the estimate \cref{29-3-19ct3} is fulfilled, and have the following result.

\begin{theorem}\label{con.rate}
	Assume that $\tilde{f} \in \mathcal{I} \left(Z\right)$ and there exists a function $w \in L^2(\Sigma)$ such that 
	\begin{align}\label{17-6-20ct1}
	\tilde{f}=F(w)+f^*,
	\end{align}
where $F(w)$ is the weak solution of the equation \cref{10-12-19ct1}. Then:
	
	(i) $\tilde{f} = f^\dag$, i.e. it is the unique solution of the identification problem $\left(\mathcal{IP}\right)$.
	
	(ii) The estimate
	\begin{align*}
	\|f-f^\dag\|^2_{L^2(\Omega_T)} \le C\left( h^3\rho^{-1} + \tau^2h^{-1}\rho^{-1} +\delta +\rho + \delta^2\rho^{-1}\right)
	\end{align*}
	holds, where  $f$ denotes the unique minimizer of $\big(\mathcal{P}_{\rho,\delta,h,\tau}\big)$.
	
	(iii) With the noise level $\delta$, the choice $h \sim \delta$, $\tau \sim \delta^{3/2}$ and $\rho \sim \delta$ leads to the convergence rate
		\begin{align*}
		\|f-f^\dag\|_{L^2(\Omega_T)} = \mathcal{O}(\sqrt{\delta}).
		\end{align*}
\end{theorem}

Assume that the coefficients appearing in the elliptic differentiable operator $\mathcal{L}$ and the domain $\Omega$ are smooth enough. Then the solution $F(w)$ to \cref{10-12-19ct1} has the regularity property $F(w) \in H^k(0,T; L^2(\Omega))$ with $\frac{\partial^{k+1} F(w)}{\partial t^{k+1}} \in L^2(0,T; H^{-1}(\Omega))$, $k \ge 0$, provided $w$ has the similar regularity (cf.\ \cite{evan,wolka}). Therefore, with an a priori smooth estimate $f^*$ the identification $f^\dag = F(w) + f^*$ and the exact $u(f^\dag)$ satisfy the inequality \cref{29-3-19ct3}. We also wish to mention that the adjoint approach in the present
paper yields the convergence rate $\mathcal{O}(\sqrt{\delta})$ for the Tikhonov regularization under the source condition \cref{17-6-20ct1}. Other source conditions may provide different convergence rates or even the optimal-order  one
$\mathcal{O}(\delta^{2/3})$, which are unfortunately still open to us.

\begin{proof}
	(i) For all  $\theta \in \mathcal{I} \left(Z\right)$ we rewrite
	\begin{align*}
	\|\theta-f^*\|^2_{L^2(\Omega_T)}- \|\tilde{f}-f^*\|^2_{L^2(\Omega_T)} 
	&= \|\theta-\tilde{f}\|^2_{L^2(\Omega_T)} + 2(\tilde{f}-f^*,\theta-\tilde{f})_{L^2(\Omega_T)} \\
	&\ge 2(\tilde{f}-f^*,\theta-\tilde{f})_{L^2(\Omega_T)}
	\end{align*}
	and so need to show that $(\tilde{f}-f^*,\theta-\tilde{f})_{L^2(\Omega_T)} =0$. In fact, by \cref{17-10-16ct2}, we have
	\begin{align*}
	&(\tilde{f}-f^*,\theta-\tilde{f})_{L^2(\Omega_T)} \\
	&~\quad = \int_0^T \int_\Omega \theta F(w)dxdt +  \int_0^T \int_{\partial\Omega} g F(w)dxdt \\
	&~\qquad - \int_0^T \int_\Omega \tilde{f} F(w)dxdt -  \int_0^T \int_{\partial\Omega} g F(w)dxdt \\
	&~\quad =\int_0^T \left\langle \frac{\partial \left( u(\theta) -u(\tilde{f})\right) }{\partial t}, F(w) \right\rangle_{\left( H^1(\Omega)^*, H^1(\Omega)\right) }dt + \int_0^T a\left( u(\theta) -u(\tilde{f}),F(w)\right) dt.
	\end{align*}
	Since $F(x,T;w)= u(x,0;\theta) -u(x,0;\tilde{f})=0$, we thus have that
	\begin{align*}
	&(\tilde{f} - f^*,\theta-\tilde{f})_{L^2(\Omega_T)}\\
	&~\quad =-\int_0^T \left\langle \frac{\partial F(w)}{\partial t}, u(\theta) -u(\tilde{f}) \right\rangle_{\left( H^1(\Omega)^*, H^1(\Omega)\right) }dt + \int_0^T a\left( F(w),u(\theta) -u(\tilde{f})\right) dt\\
	&~\quad =\int_\Sigma w(u(\theta) -u(\tilde{f}))dxdt =0,
	\end{align*}
	by $\tilde{f}, \theta \in \mathcal{I} \left(Z\right)$.
	
	(ii) By the optimality of $f$, we get that
	\begin{align}\label{29-3-19ct4}
	&\sum_{n=1}^M \int_{t^{n-1}}^{t^n}\|U^n(f)-z_\delta\|^2_{L^2(\Gamma)} dt+ \rho\|f-f^*\|^2_{L^2(\Omega_T)}\notag\\
	&~\quad \le \sum_{n=1}^M \int_{t^{n-1}}^{t^n}\|U^n(f^\dag)-z_\delta\|^2_{L^2(\Gamma)} dt + \rho\|f^\dag-f^*\|^2_{L^2(\Omega_T)}
	\end{align}
	with
	\begin{align}\label{29-3-19ct5}
	\sum_{n=1}^{M} \int_{t^{n-1}}^{t^n}\|U^n(f^\dag)-z_{\delta}\|^2_{L^2(\Gamma)} dt
	&= \sum_{n=1}^{M} \int_{t^{n-1}}^{t^n}\|U^n(f^\dag)- u(f^\dag) + Z -z_{\delta}\|^2_{L^2(\Gamma)} dt \notag\\
	&\le C\left( \sum_{n=1}^{M} \int_{t^{n-1}}^{t^n}\|U^n(f^\dag)-u(f^\dag)\|^2_{L^2(\Gamma)}dt + \delta^2\right).
	\end{align}
	Furthermore, by \cref{17-10-16ct5***}, we get that
	\begin{align}\label{29-3-19ct6}
	\sum_{n=1}^{M} \int_{t^{n-1}}^{t^n}\|U^n(f^\dag)-u(f^\dag)\|^2_{L^2(\Gamma)}dt 
	& \le Ch^{-1} \sum_{n=1}^{M} \int_{t^{n-1}}^{t^n}\|U^n(f^\dag)-u(f^\dag)\|^2_{L^2(\Omega)}dt\notag\\
	&~\quad  + Ch \sum_{n=1}^{M} \int_{t^{n-1}}^{t^n}\|U^n(f^\dag)-u(f^\dag)\|^2_{H^1(\Omega)}dt.
	\end{align}
	Due to \cref{29-3-19ct3}, the first term in the right hand side of \cref{29-3-19ct6} is bounded by
	\begin{align}\label{2-4-19ct1}
	&h^{-1} \sum_{n=1}^{M} \int_{t^{n-1}}^{t^n}\|U^n(f^\dag)-u(f^\dag)\|^2_{L^2(\Omega)}dt \notag\\
	&~\quad \le Ch^{-1} \Bigg( \sum_{n=1}^{M} \int_{t^{n-1}}^{t^n}\|U^n(f^\dag)-u^n(f^\dag)\|^2_{L^2(\Omega)}dt \notag\\
	&~\qquad\qquad\qquad\quad+ \sum_{n=1}^{M} \int_{t^{n-1}}^{t^n}\|u^n(f^\dag) - u(f^\dag)\|^2_{L^2(\Omega)}dt\Bigg) \notag\\
	&~\quad \le Ch^{-1} \left( h^4 + \tau^4 + \sum_{n=1}^{M} \int_{t^{n-1}}^{t^n}\|u^n(f^\dag) - u(f^\dag)\|^2_{L^2(\Omega)}dt\right)\notag\\
	&~\quad \le Ch^{-1} \left( h^4 + \tau^4 + \tau^2 \|u_{t}(f^\dag)\|^2_{L^\infty(0,T; L^2(\Omega))}\right)\notag\\
	&~\quad \le C(h^3 + \tau^4h^{-1} +\tau^2h^{-1}),
	\end{align}
	where we used the estimate
		\begin{align*}
		\sum_{n=1}^{M} \int_{t^{n-1}}^{t^n}\|u^n(f^\dag) - u(f^\dag)\|^2_{L^2(\Omega)}dt 
		&= \sum_{n=1}^{M} \int_{t^{n-1}}^{t^n}\|u_t(s) (t^n - t^{n-1})\|^2_{L^2(\Omega)}dt \\
		&\le \tau^2 \|u_{t}(f^\dag)\|^2_{L^\infty(0,T; L^2(\Omega))}.
		\end{align*}
	Likewise, we have that
	\begin{align}\label{2-4-19ct2}
	h \sum_{n=1}^{M} \int_{t^{n-1}}^{t^n}\|U^n(f^\dag)-u(f^\dag)\|^2_{H^1(\Omega)}dt
	&\le Ch \left( h^2 + \tau^4h^{-2} +\tau^2\|u_{t}(f^\dag)\|^2_{L^\infty(0,T; H^1(\Omega))}\right) \notag\\
	&\le C(h^3 +\tau^4h^{-1} +\tau^2h).
	\end{align}
	We deduce from \cref{29-3-19ct4}--\cref{2-4-19ct2} that
	\begin{align}\label{20-3-19ct1}
	&\sum_{n=1}^M \int_{t^{n-1}}^{t^n}\|U^n(f)- z_\delta\|^2_{L^2(\Gamma)}dt + \rho\|f-f^\dag\|^2_{L^2(\Omega_T)}\notag\\
	&~\quad \le C\left( h^3  +\tau^2h^{-1}  +\delta^2\right)   + 2\rho(f^\dag-f^*,f^\dag-f)_{L^2(\Omega_T)}.
	\end{align}
	Since $f^\dag - f^* = F(w)$, we have
	\begin{align}\label{29-3-19ct8}
	&(f^\dag-f^*,f^\dag-f)_{L^2(\Omega_T)} \notag\\
	&~\quad = \int_\Sigma w(u(f^\dag) -u(f))dxdt \notag\\
	&~\quad = \int_\Sigma w(u(f^\dag) - z_\delta)dxdt + \int_\Gamma \left( \sum_{n=1}^{M} \int_{t^{n-1}}^{t^n} w(z_\delta - U^n(f))dt\right) dx \notag\\
	&~\qquad + \int_\Gamma \left( \sum_{n=1}^{M} \int_{t^{n-1}}^{t^n} w( U^n(f) - u(f))dt\right) dx.
	\end{align}
	We bound for each term in the right hand side of the above equation. First, we get
	\begin{align}\label{29-3-19ct9}
	\int_\Sigma w(u(f^\dag) - z_\delta)dxdt
	&\le \|w\|_{L^2(\Sigma)} \|u(f^\dag) - z_\delta\|_{L^2(\Sigma)} \le \|w\|_{L^2(\Sigma)}\delta
	\end{align}
	and
	\begin{align}\label{29-3-19ct10}
	&\int_\Gamma \left( \sum_{n=1}^{M} \int_{t^{n-1}}^{t^n} w(z_\delta - U^n(f))dt\right) dx \notag\\
	&~\quad \le \int_\Gamma \left( \sum_{n=1}^{M} \left( \int_{t^{n-1}}^{t^n} w^2dt\right)^{1/2} \left( \int_{t^{n-1}}^{t^n} (z_\delta - U^n(f))^2dt\right)^{1/2} \right) dx \notag\\
	&~\quad \le \int_\Gamma \left(  \left( \sum_{n=1}^{M}\int_{t^{n-1}}^{t^n} w^2dt\right)^{1/2} \left( \sum_{n=1}^{M}\int_{t^{n-1}}^{t^n} (z_\delta - U^n(f))^2dt\right)^{1/2} \right) dx \notag\\
	&~\quad \le  \left( \int_\Gamma \left(\sum_{n=1}^{M}\int_{t^{n-1}}^{t^n} w^2dt \right)dx\right)^{1/2} \left( \int_\Gamma \left( \sum_{n=1}^{M}\int_{t^{n-1}}^{t^n} (z_\delta - U^n(f))^2dt\right)dx\right)^{1/2} \notag\\
	&~\quad \le \rho\|w\|^2_{L^2(\Sigma)} + \dfrac{1}{4\rho}\sum_{n=1}^M \int_{t^{n-1}}^{t^n}\|U^n(f)-z_\delta\|^2_{L^2(\Gamma)}dt.
	\end{align}
	Further, as the above estimates \cref{29-3-19ct6}--\cref{2-4-19ct2}, we get
	\begin{align}\label{29-3-19ct11}
	\int_\Gamma \left( \sum_{n=1}^{M} \int_{t^{n-1}}^{t^n} w( U^n(f) - u(f))dt\right) dx \le C\left( h^3 +\tau^2h^{-1}\right).
	\end{align}
	It follows from \cref{29-3-19ct8}--\cref{29-3-19ct11} that
	\begin{align}\label{29-3-19ct12}
	&2\rho(f^\dag-f^*,f^\dag-f)_{L^2(\Omega_T)}\notag\\
	&~\quad \le C\left( h^3 +\tau^2h^{-1} + \rho\delta +\rho^2\right) + \dfrac{1}{2}\sum_{n=1}^M \int_{t^{n-1}}^{t^n}\|U^n(f)-z_\delta\|^2_{L^2(\Gamma)}dt.
	\end{align}
	We therefore conclude from \cref{20-3-19ct1} and \cref{29-3-19ct12}
	\begin{align*}
	\dfrac{1}{2}\sum_{n=1}^M \int_{t^{n-1}}^{t^n}\|U^n(f)- z_\delta\|^2_{L^2(\Gamma)} + \rho\|f-f^\dag\|^2_{L^2(\Omega_T)} \le C\left( h^3 +\tau^2h^{-1} +\rho\delta+\rho^2 +\delta^2\right),
	\end{align*}
	which completes the proof.
\end{proof}

\section{Numerical examples} \label{NumerEx}
In this section we employ the conjugate gradient (CG) method (cf.\ \cite{Kelley99}) for reaching the unique minimizer of the finite dimensional minimization problem $(\mathcal{P}_{\rho,\delta,h,\tau})$. Starting with an initial guess $f_0$, the iterative process is given by
\begin{equation}\label{CGupdatestep}
f_{k+1} = f_k + \alpha_kd_k,\quad k = 0,1,\ldots
\end{equation}
where $d_k$ is the update direction and $\alpha_k$ is the update step size until it meets a stopping criterion of the form
\begin{equation}\label{CGstoppingcond}
\|\nabla J_{\rho,\delta,h,\tau}(f_{k})\| \leq \tau_a + \tau_r\|\nabla J_{\rho,\delta,h,\tau}(f_{0})\|.
\end{equation}
The update direction is computed as
\begin{equation}\label{CGupdatedirect}
d_k = \begin{cases}
-\nabla J_{\rho,\delta,h,\tau}(f_{k}) &\quad \mbox{ if } \quad k = 0,\\
-\nabla J_{\rho,\delta,h,\tau}(f_{k})+\beta_kd_{k-1} &\quad \mbox{ if } \quad k > 0,
\end{cases}
\end{equation}
where coefficient $\beta_k$ is computed using the Polak-Ribi\`ere formula
\begin{equation}\label{CGPolakR}
\beta_k = \dfrac{\left( \nabla J_{\rho,\delta,h,\tau}(f_{k}),\nabla J_{\rho,\delta,h,\tau}(f_{k})-\nabla J_{\rho,\delta,h,\tau}(f_{k-1})\right)_{{L^2(\Omega_T)}^d}}{\|\nabla J_{\rho,\delta,h,\tau}(f_{k-1})\|^2_{L^2(\Omega_T)}}.
\end{equation}
To compute $\alpha_k$, we consider the quadratic minimization problem
\begin{equation}\label{CGquadoptimprob}
\arg\min_{\alpha \geq 0}J_{\rho,\delta,h,\tau}(f_k + \alpha d_k).
\end{equation}
It is shown that the unique minimizer of \cref{CGquadoptimprob} is given by
\begin{equation}\label{CGquadoptimminimizer}
\alpha_k = -\dfrac{\sum_{n=1}^M \int_{t^{n-1}}^{t^n} \left ( U^n_0(d_k) , U^n(f_k) -z_\delta \right)_{L^2(\Gamma)}dt + \rho ( d_k,f_k-f^* )_{L^2(\Omega_T)}}{\sum_{n=1}^M \tau \|U^n_0(d_k)\|^2_{L^2(\Gamma)} + \rho\|d_k\|^2_{L^2(\Omega_T)}},
\end{equation}
where for each $n\in I$, $U^n_0(d_k) \in \mathcal{V}^1_h$ is the unique solution of the equation
\begin{align*}
\left( \partial U^n_0(d_k), \varphi_h \right)_{L^2(\Omega)} + a^n(\overline{\partial}U^n_0(d_k),\varphi_h) &= \left(\overline{d_k}^n,\varphi_h \right)_{L^2(\Omega)} \quad \forall \varphi_h\in \mathcal{V}^1_h,\\
U^0_0(f) &=0,
\end{align*}
in which $\left( U^n(f_k)\right)_{n=1}^M$ is defined by the system \cref{18-1-17ct5*}. 
Practical steps can be summarized in \Cref{Alg:CG4IP}.

\begin{algorithm}[H]
	\caption{CG method for minimizing the cost functional $\left(\mathcal{P}_{\rho,\delta,h,\tau}\right)$}
	\label{Alg:CG4IP}
	\begin{algorithmic}[1]
		\REQUIRE  $f_0, \tau_r, \tau_a, k_{\max}$
		\ENSURE An approximate of the minimizer to $\left(\mathcal{P}_{\rho,\delta,h,\tau}\right)$
		\STATE  Calculate $\nabla J_{\rho,\delta,h,\tau}(f_{0})$ and set $d_0 =-\nabla J_{\rho,\delta,h,\tau}(f_{0})$.
		\STATE  Compute $\alpha_0$ using \cref{CGquadoptimminimizer}.
		\STATE  Update $f_1 = f_0 + \alpha_0d_0$.
		\STATE Set $k=1$ and compute $\nabla J_{\rho,\delta,h,\tau}(f_{k})$.
		\WHILE{$k \le k_{\max}$ and \cref{CGstoppingcond} does not hold} \STATE Compute $\beta_k$, using \cref{CGPolakR}.
		\STATE Update the direction, using \cref{CGupdatedirect}.
		\STATE Compute $\alpha_k$, using \cref{CGquadoptimminimizer}.
		\STATE Update the minimizer, using \cref{CGupdatestep}.
		\STATE Set $k=k+1$ and compute $\nabla J_{\rho,\delta,h,\tau}(f_{k})$
		\ENDWHILE
	\end{algorithmic}
\end{algorithm}

For numerical implementations below, a common setting is used, unless stated otherwise. We consider a spatial two-dimensional case of \cref{ct1}, where $\Omega_T = (-1,1)\times (-1,1)  \times (0,1]$ and
\begin{align}\label{17-6-20ct2}
A=\left[\begin{matrix}3&1\\ 1 & 2\end{matrix}\right], b = 1, \sigma = 1.
\end{align}
We denote by $f_{mean}$ the mean of the exact source $f$. For initial iterate and predicted source, we respectively choose $f_0 = 0$ and $f^* = f + 0.2(f-f_{mean})$. To mimic the noise measurement, we set $z_{\delta} = u(f)_{|\Sigma}+c_w  z_{rand}$ where $z_{rand}$ is produced by MATLAB command `rand' and $c_w$ is chosen such that the measurement error $\|u(f)-z_\delta\|_{L^2(\Sigma)}$ is equal to the prescribed noise level $\delta$.

Except for the first example where we choose different mesh sizes to show the experimental order of convergence (EOC), we will take $h = 0.05$ and generate the mesh using MATLAB function 'generateMesh'. Choice for other constant parameters conceptually follows from  \Cref{stability2}
\begin{equation*}
\tau = 0.25\cdot h,\quad \rho = 0.01\cdot h,\quad \delta = 0.5 \cdot h^{2}.
\end{equation*}

\subsection{Time dependent source $f(t)$}\label{Subsec:timdependEx}
\begin{figure}[H]
\begin{center}
	\begin{minipage}{.33\textwidth}
		\includegraphics[width=\textwidth]{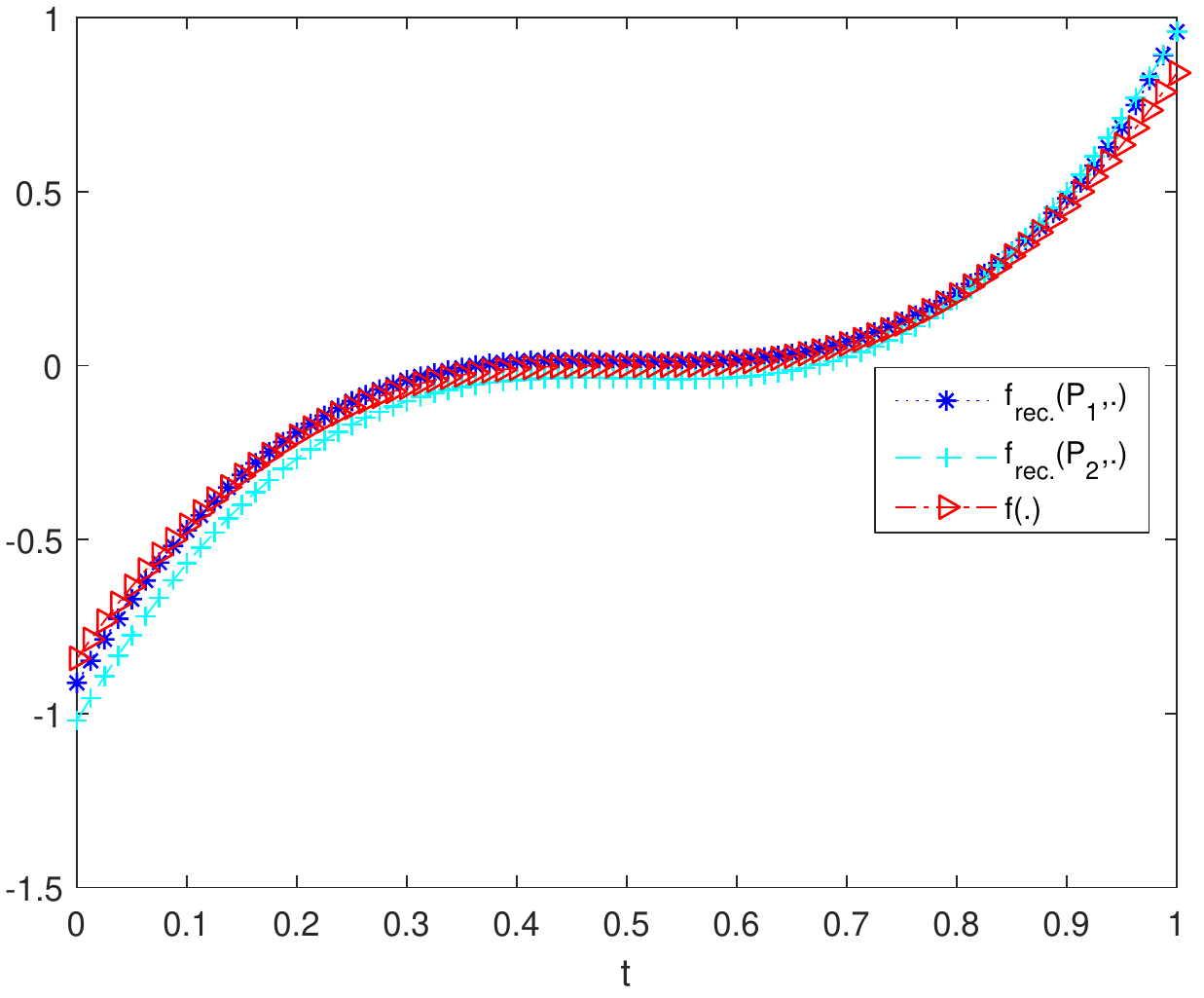}\\
		\hspace*{.45\textwidth} {\footnotesize (a)}
	\end{minipage}
	\begin{minipage}{.33\textwidth}
		\includegraphics[width=\textwidth]{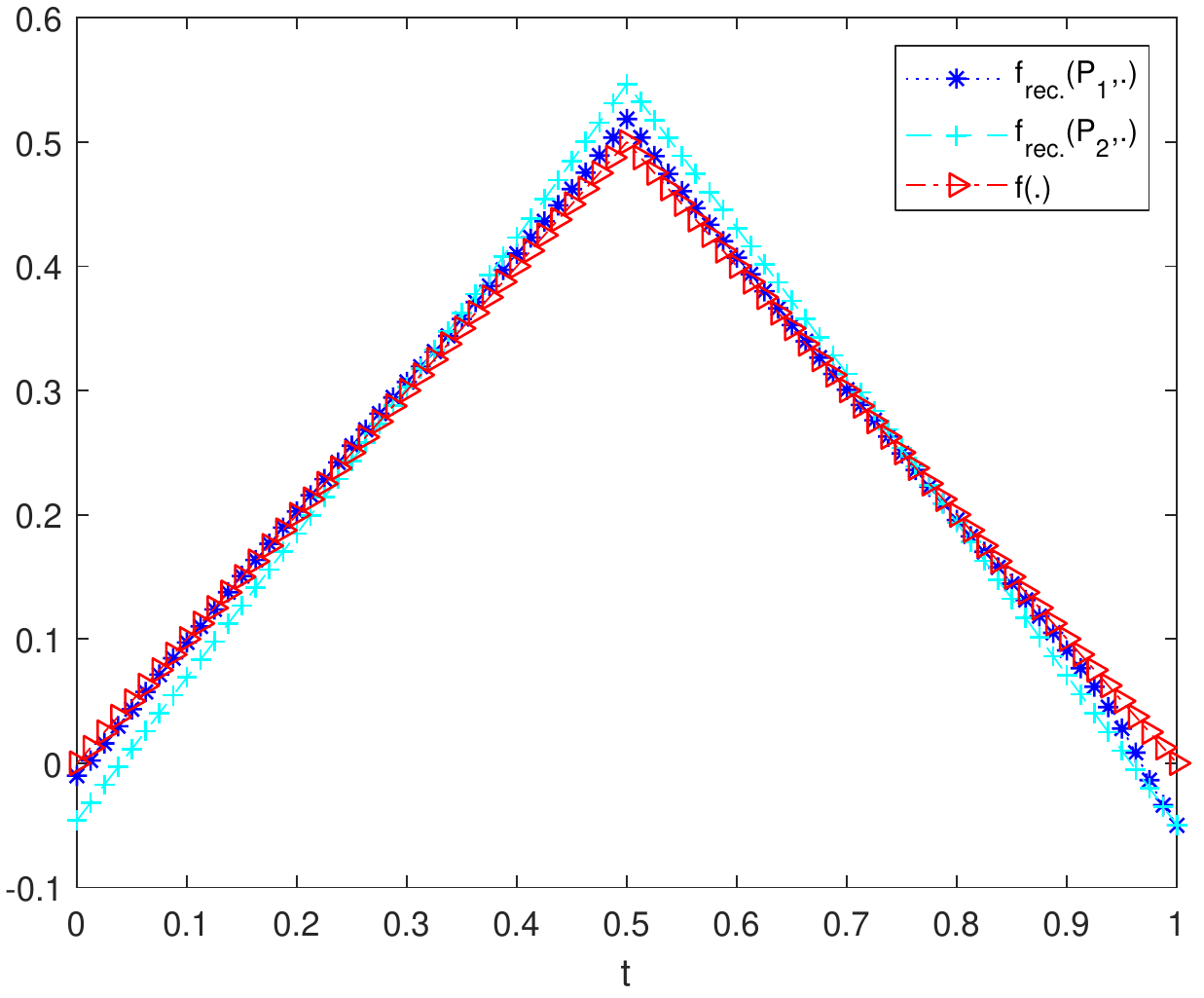} \\
		\hspace*{.45\textwidth} {\footnotesize (b)}
	\end{minipage}
	\begin{minipage}{.33\textwidth}
		\includegraphics[width=\textwidth]{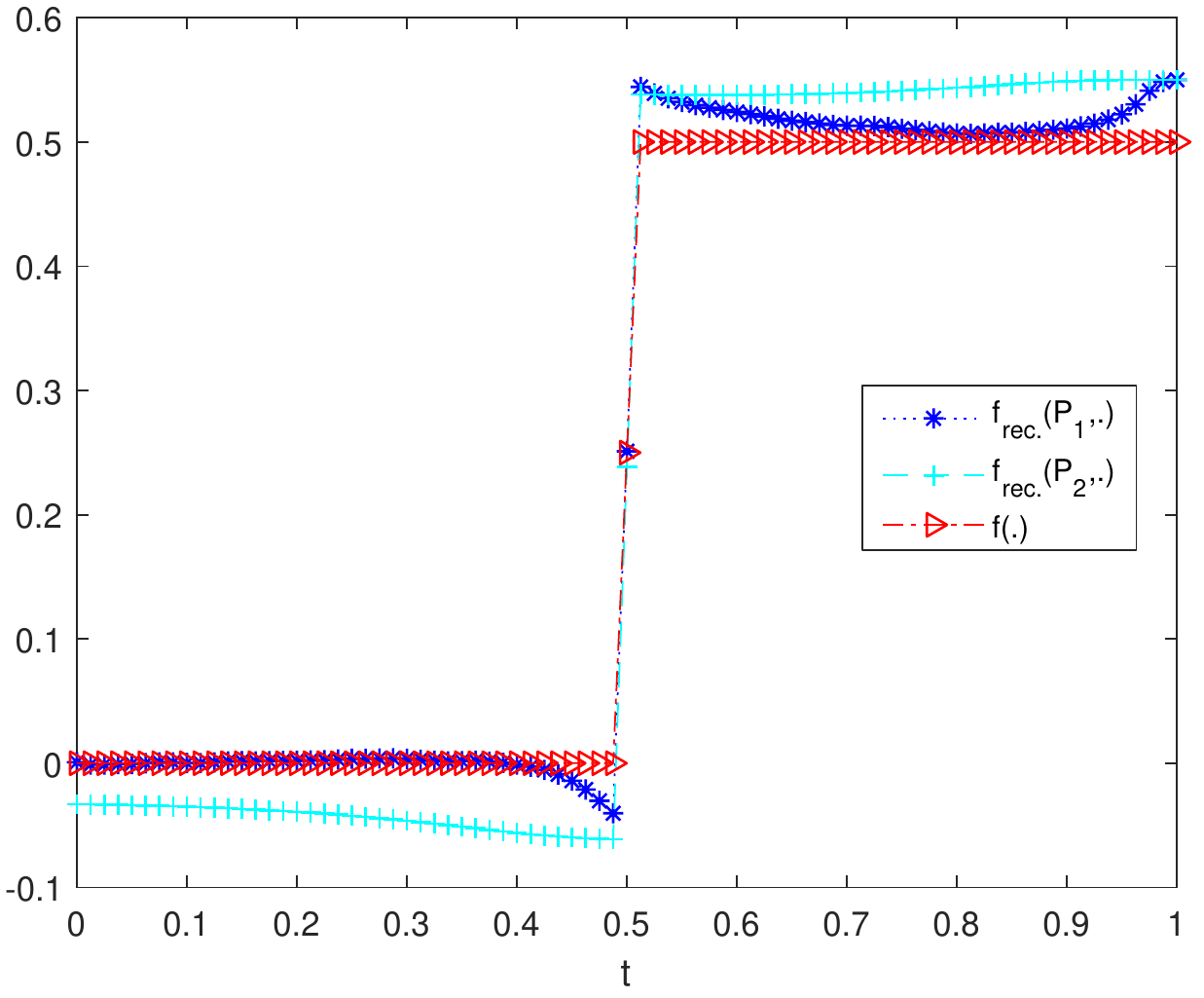} \\
		\hspace*{.45\textwidth} {\footnotesize (c)}
	\end{minipage}
	\end{center}
	\caption{Time dependent source: comparisons of the recovered source $f_{rec.}$ and exact source at $P_1(-0.1,-0.5)$ and $P_2(0.5,0.6)$ in different cases: $f(t) =(2t-1)^2 \sin(2t-1)$ (a), $f(t) = 0.5 - |0.5-t|$ (b), $f(t) = 0.5H(t-0.5)$ (c), where $H$ is the Heaviside step function.}
	\label{Fig:Mod3_comp_source}
\end{figure}

In this subsection, we specify the right hand side $f(t)$ independent of the spatial variables, the initial value $u(\cdot,\cdot,0) = 0.4$, and boundary condition $g = 0.4$. Then the state is numerically computed by solving the corresponding forward problem. We will refer to this source function and the resulting numerical state as the exact source and the exact state, respectively. To generate the noise data $z_\delta$, we simulate the forward problem with the aforementioned data but on a finer grid with $h_{fine}=0.025$ to prevent the so called inverse crime (cf.\ \cite{Wirg04}). The derived solution is then interpolated on the grid with mesh size $h$ and perturbed.
The observation boundary part is $\Sigma = [-1,1]\times \{-1\} \times (0,1]$ and the result is checked at $P_1$ and $P_2$, defined as the closest nodal point to $(-0.1,-0.5)$ and $(0.5,0.6)$, respectively. In \Cref{Fig:Mod3_comp_source} one can see that at both $P_1$ and $P_2$, the recovered source matches the exact source quite well in various cases.

\begin{figure}[H]
\begin{center}
	\begin{minipage}{.33\textwidth}
		\includegraphics[width=\textwidth]{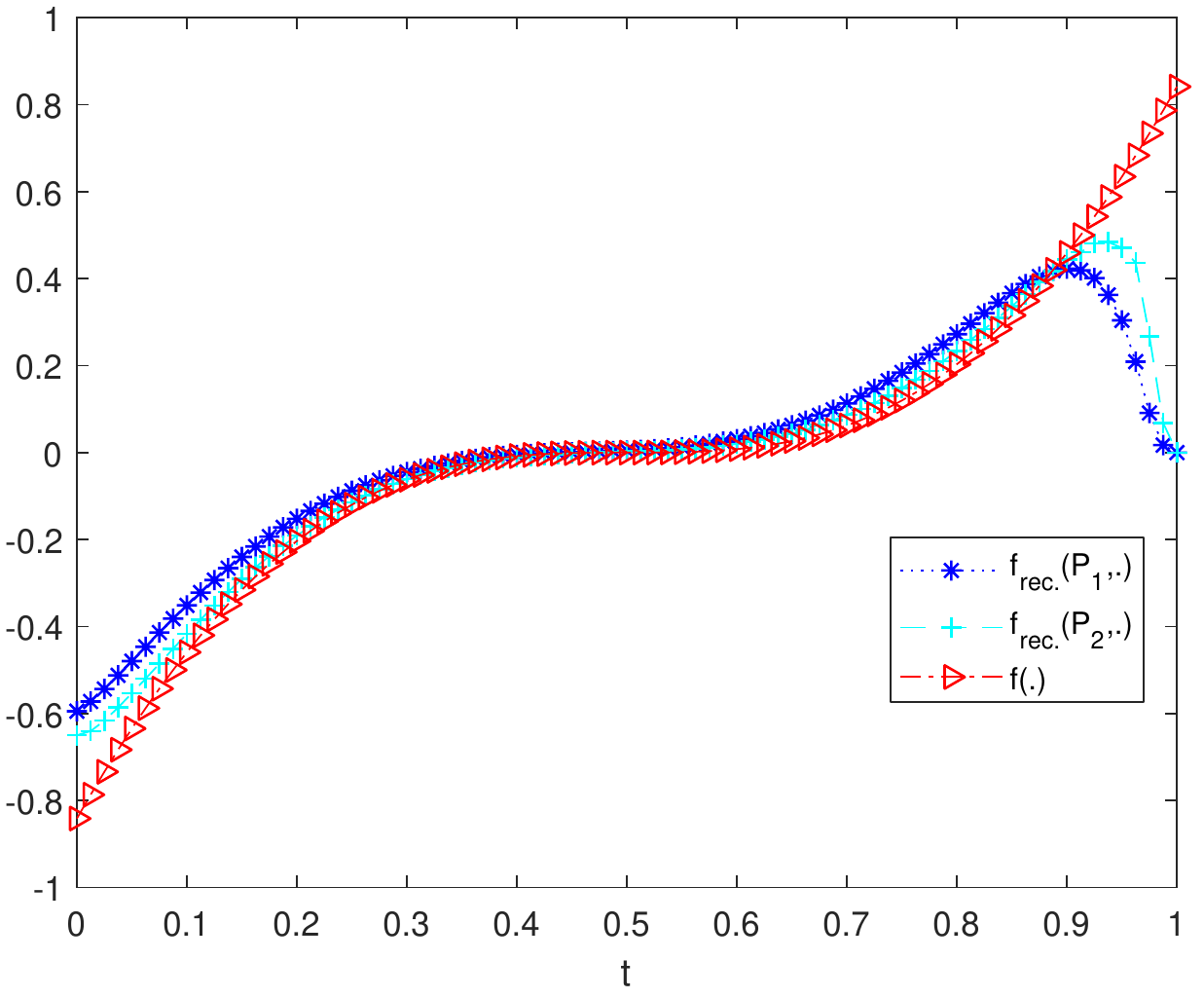}\\
		\hspace*{.45\textwidth} {\footnotesize (a)}
	\end{minipage}
	\begin{minipage}{.33\textwidth}
		\includegraphics[width=\textwidth]{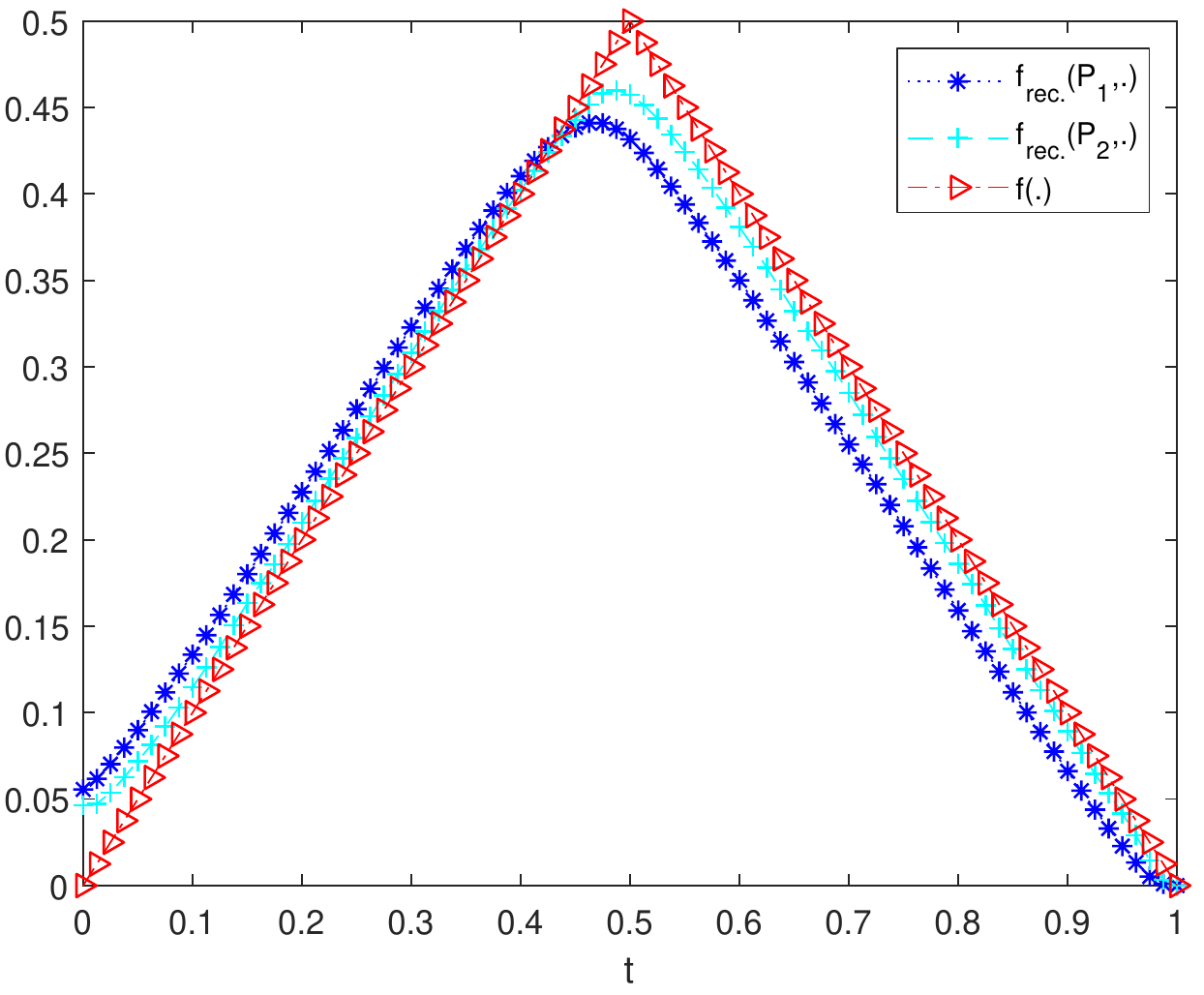} \\
		\hspace*{.45\textwidth} {\footnotesize (b)}
	\end{minipage}
	\begin{minipage}{.33\textwidth}
		\includegraphics[width=\textwidth]{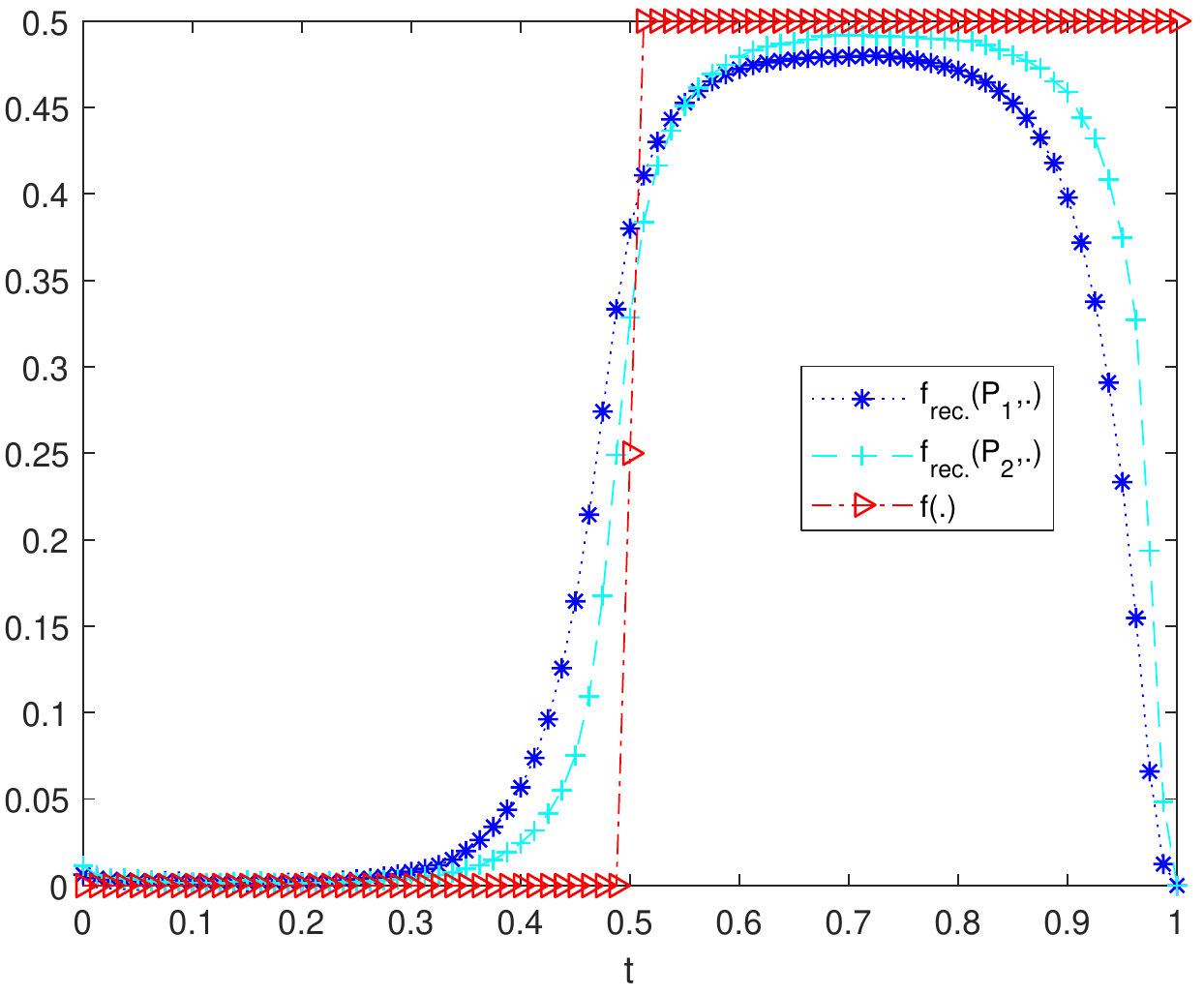} \\
		\hspace*{.45\textwidth} {\footnotesize (c)}
	\end{minipage}
	\end{center}
	\caption{Time dependent source with $f^* = 0$: comparisons of the recovered source and exact source at $P_1(-0.1,-0.5)$ and $P_2(0.5,0.6)$ in different cases: $f(t) =(2t-1)^2 \sin(2t-1)$ (a), $f(t) = 0.5 - |0.5-t|$ (b), $f(t) = 0.5H(t-0.5)$ (c), where $H$ is the Heaviside step function.}
	\label{Fig:Mod3_comp_source2}
\end{figure}

We would like to mention that for this special case, the considered inverse problem has a unique solution (cf.\ \cite{HaPe13,Slodika2015}). That is, wherever we choose an a priori estimate $f^*$, the algorithm should always yield a good approximation of the exact source. In practice, a poorly predicted source may however reduce the quality of the result. In \Cref{Fig:Mod3_comp_source2} we perform the same test as in \Cref{Fig:Mod3_comp_source} but with $f^* = 0$ for all cases. One can see that the obtaining numerical solutions are not so good, compared with the previous simulations (cf.\ \Cref{Fig:Mod3_comp_source}). 

\subsection{Space-dependent source $f(x,y)$}\label{Subsec:DiscontEx}
This implementation is constructed similar to that of \Cref{Subsec:timdependEx}, where the source function is independent of the time variable and chosen as
\begin{equation*}
f(x,y,t)=\left\{\begin{array}{ll}
0.5 &\ \mbox{if} \quad (x,y)\in B(O,0.5) :=\{(x,y):x^2 + y^2 <0.5^2\},\\
0 &\ \mbox{otherwise}
\end{array}\right.
\end{equation*}
for all $t$. We note that if the source term depends on the spatial variables, certain space observations are required to guarantee the uniqueness of the identification problem, for example, the additional  final data measurement. For this subject we mention to
\cite{FaLe06,Ha07,Ha12,JoLe07,Ru80} and references given there for detailed discussions.

\begin{figure}[H]
\begin{center}
	\begin{minipage}{.24\textwidth}
		\includegraphics[width=\textwidth]{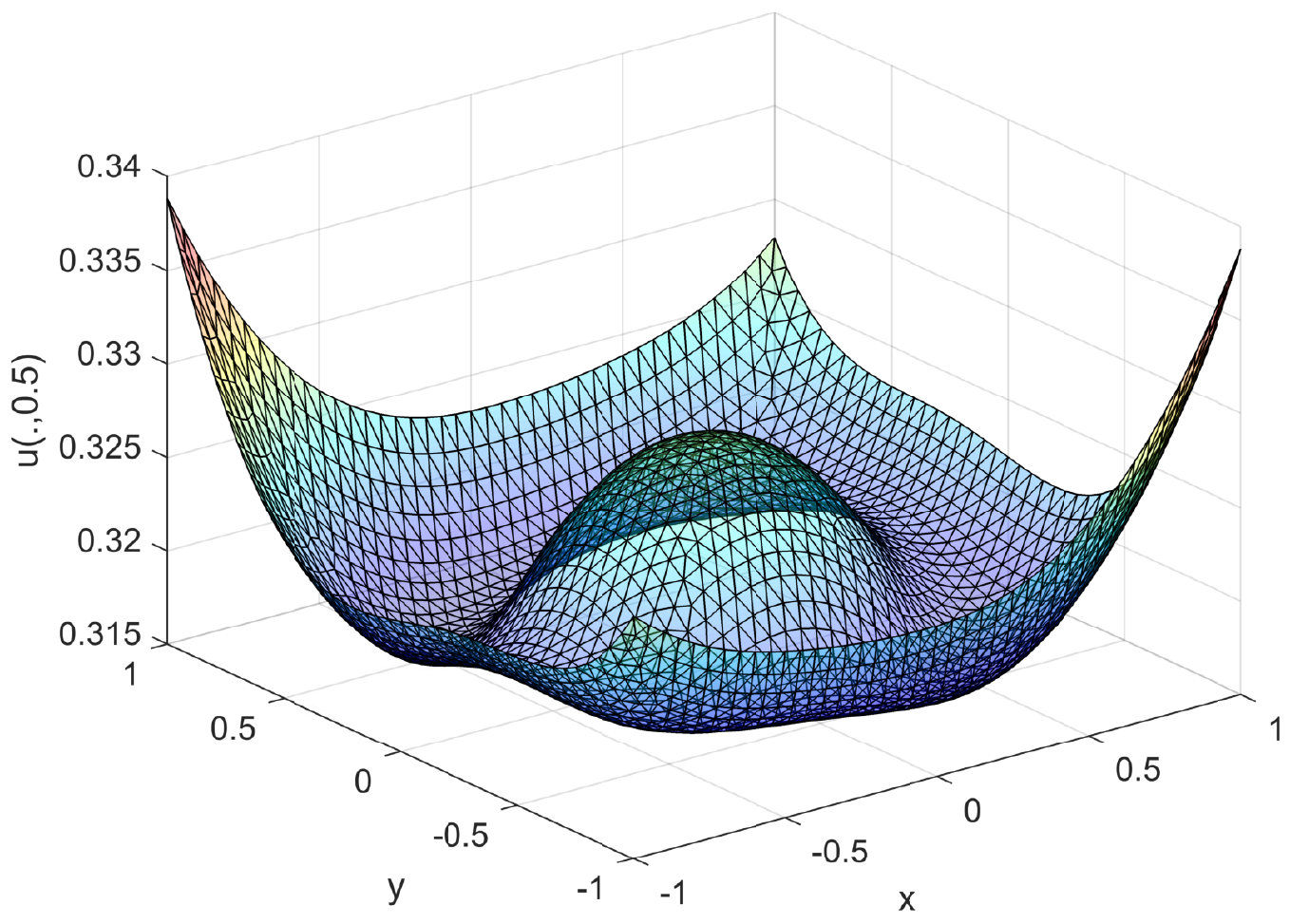}\\
		\hspace*{.4\textwidth} {\footnotesize (a)}
	\end{minipage}
	\begin{minipage}{.24\textwidth}
		\includegraphics[width=\textwidth]{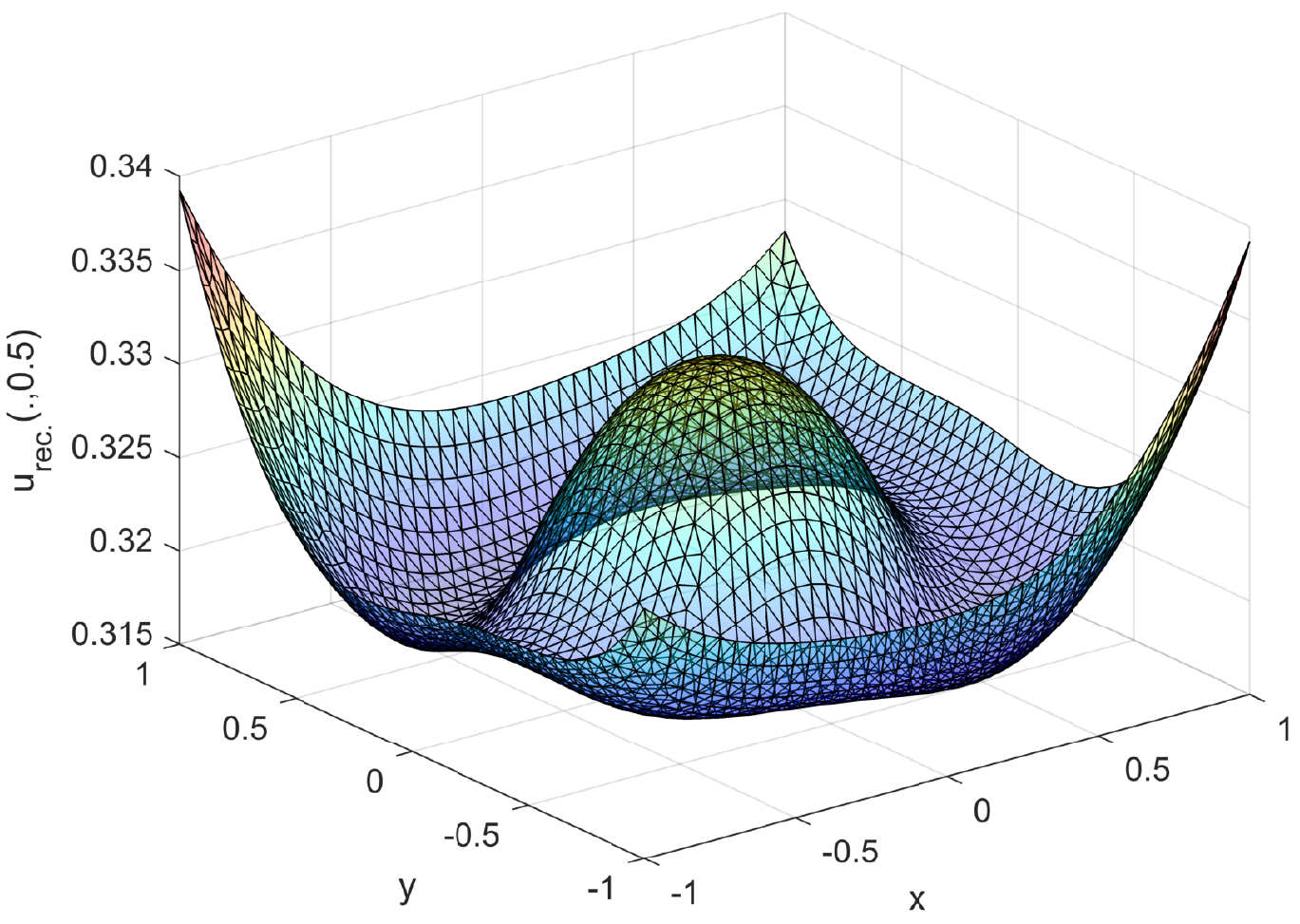}\\
		\hspace*{.4\textwidth} {\footnotesize (b)}
	\end{minipage}
		\begin{minipage}{.24\textwidth}
			\includegraphics[width=\textwidth]{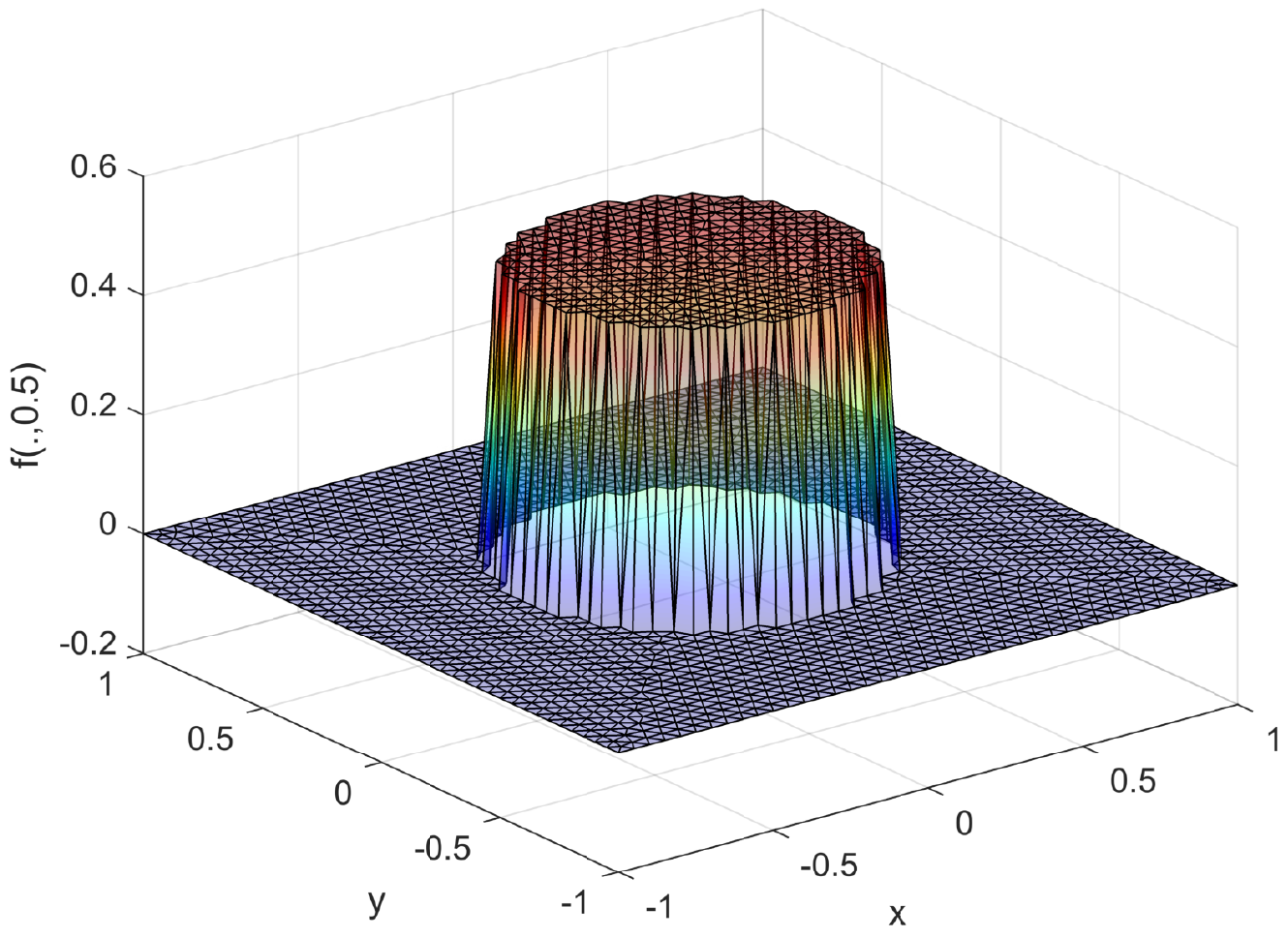}\\
			\hspace*{.4\textwidth} {\footnotesize (c)}
		\end{minipage}
		\begin{minipage}{.24\textwidth}
			\includegraphics[width=\textwidth]{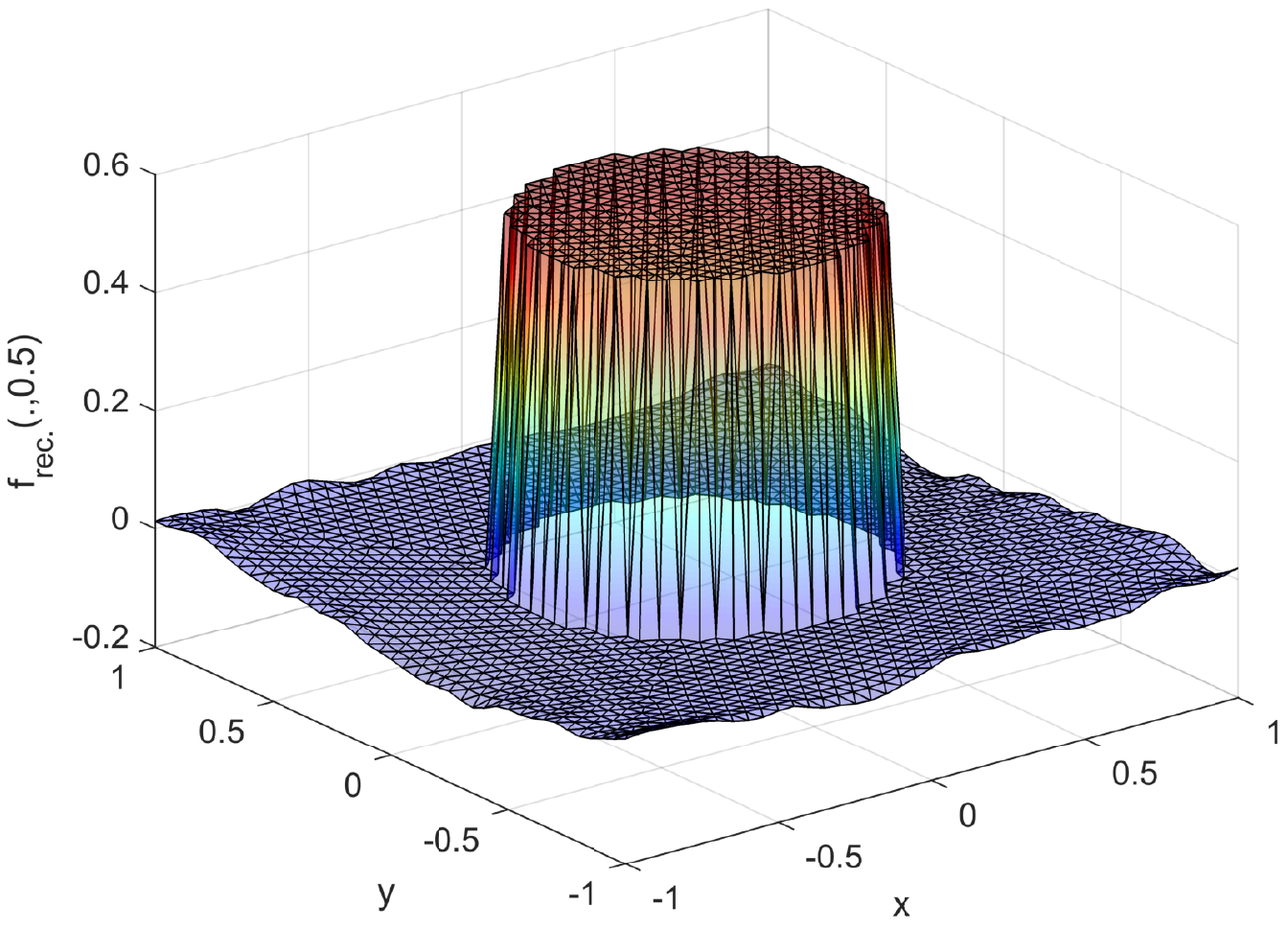}\\
			\hspace*{.4\textwidth} {\footnotesize (d)}
		\end{minipage}
		\end{center}
	\caption{Space-dependent source at instant $t=0.5$: the exact state (a), recovered state (b),  exact source (c) and recovered source (d).}
	\label{Fig:Mod2 comp_u_f}
\end{figure}

We assume that the measurement $z_\delta$ is available on the whole boundary $\Sigma = \partial \Omega \times [0,T]$.
In \Cref{Fig:Mod2 comp_u_f} we compare the recovered state/source with exact state/source with respect to the spatial variables at the time level $t=0.5$. 

We also present the test on different mesh sizes and the numerical result is shown in \Cref{Tab:Mod2MeshRefine} with its corresponding EOC being given in \Cref{Tab:Mod2EOC}.

 \begin{table}[H]
 \caption{Space-dependent source: Refinement level $l$, mesh size $h$, measurement error $\delta$, regularization parameter $\rho$, and errors}
 \begin{center}
 \begin{tabular}{ c l l l l l l}\hline
 $l$ & $h$ & $\delta$ & $\rho$& $\|u-u_{rec.}\|_{L^2(\Omega_T)}$& $\|u-u_{rec.}\|_{L^2(\Sigma)}$ &$\|f-f_{rec.}\|_{L^2(\Omega_T)}$\\  \hline
1 & 0.8 & 0.32& 0.008&0.2160& 0.2762 &1.1199\\
2 & 0.4 & 0.08& 0.004&0.0534& 0.0727 &0.3262\\
3 & 0.2 & 0.02& 0.002&0.0132& 0.0190 &0.1083\\
4 & 0.1 & 0.005&0.001&0.0029& 0.0049 &0.0556\\ \hline
 \end{tabular}
 \end{center}
 \label{Tab:Mod2MeshRefine}
 \end{table}
 
  \begin{table}[H]
  \caption{Space-dependent source: EOC}
  \begin{center}
  \begin{tabular}{ c l l l}\hline
  $l$ &  $\|u-u_{rec.}\|_{L^2(\Omega_T)}$ &$\|u-u_{rec.}\|_{L^2(\Sigma)}$ &$\|f-f_{rec.}\|_{L^2(\Omega_T)}$\\  \hline
 1 & -- & -- & -- \\
 2 & 2.0161& 1.9257 &1.7795 \\
 3 & 2.0163 &1.9360 & 1.5907\\
 4 & 2.1864&1.9551 & 0.9619\\
  Mean of EOC &2.0729 & 1.9389 & 1.4440 \\ \hline
  \end{tabular}
  \end{center}
  \label{Tab:Mod2EOC}
  \end{table}
The convergence history given in \Cref{Tab:Mod2MeshRefine} and  \Cref{Tab:Mod2EOC} shows that
the algorithm performs well for our identification problem.

\subsection{General source $f(x,y,t)$}\label{Subsec:smoothEx}
\begin{figure}[H]
\begin{center}
	\begin{minipage}{.33\textwidth}
		\includegraphics[width=\textwidth]{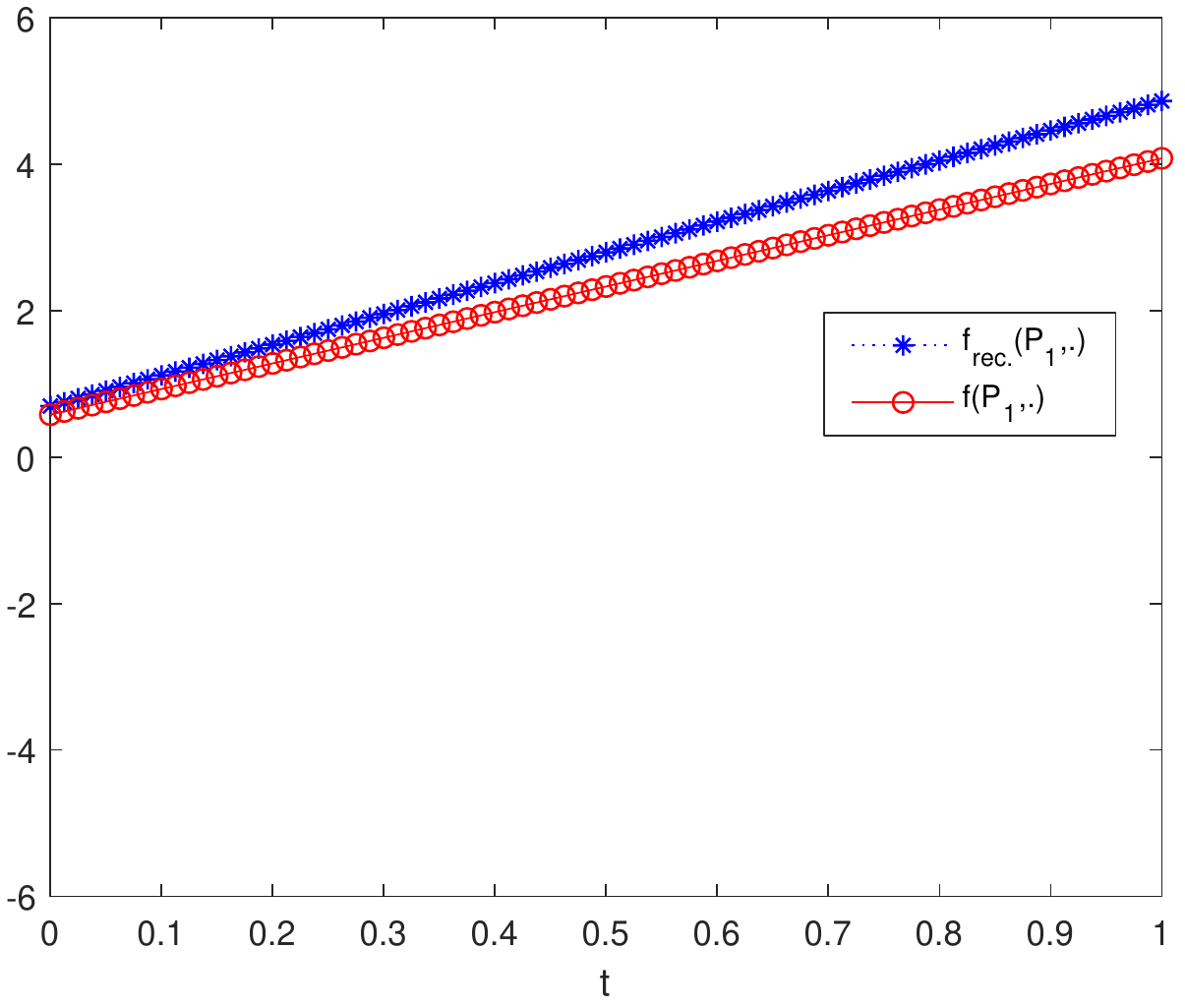}\\
		\hspace*{.45\textwidth} {\footnotesize (a)}
	\end{minipage}
	\begin{minipage}{.33\textwidth}
		\includegraphics[width=\textwidth]{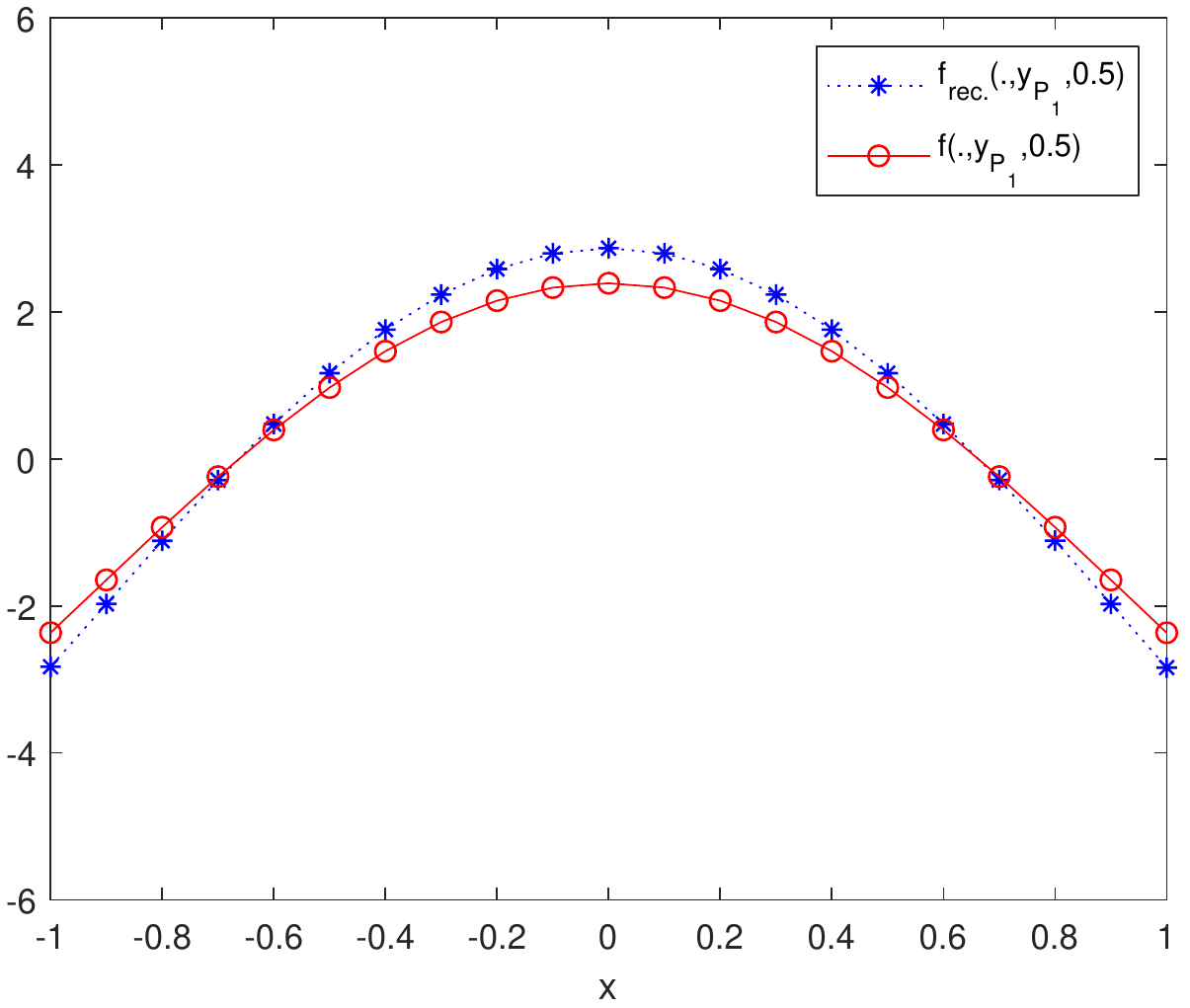} \\
		\hspace*{.45\textwidth} {\footnotesize (b)}
	\end{minipage}
	\begin{minipage}{.33\textwidth}
		\includegraphics[width=\textwidth]{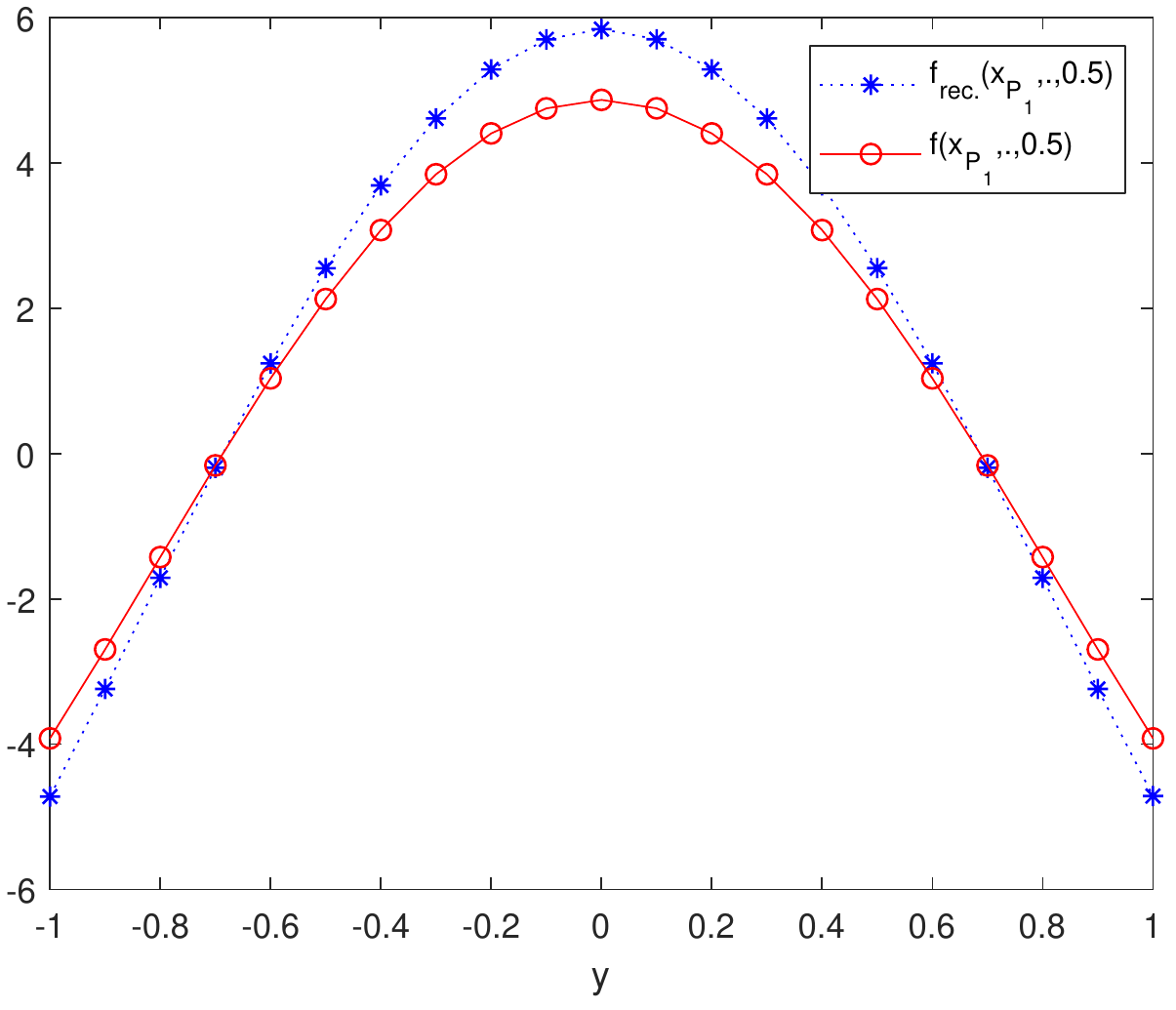} \\
		\hspace*{.45\textwidth} {\footnotesize (c)}
	\end{minipage}
	\end{center}
	\vspace*{-0.3cm}
	\caption{General source: the recovered source and exact source at $P_1$ along $t$ (a), $x$ (b), and $y$ (c).}
	\label{Fig:Mod1 comp_P1}
\end{figure}

\begin{figure}[H]
	\begin{center}
		\begin{minipage}{.33\textwidth}
			\includegraphics[width=\textwidth]{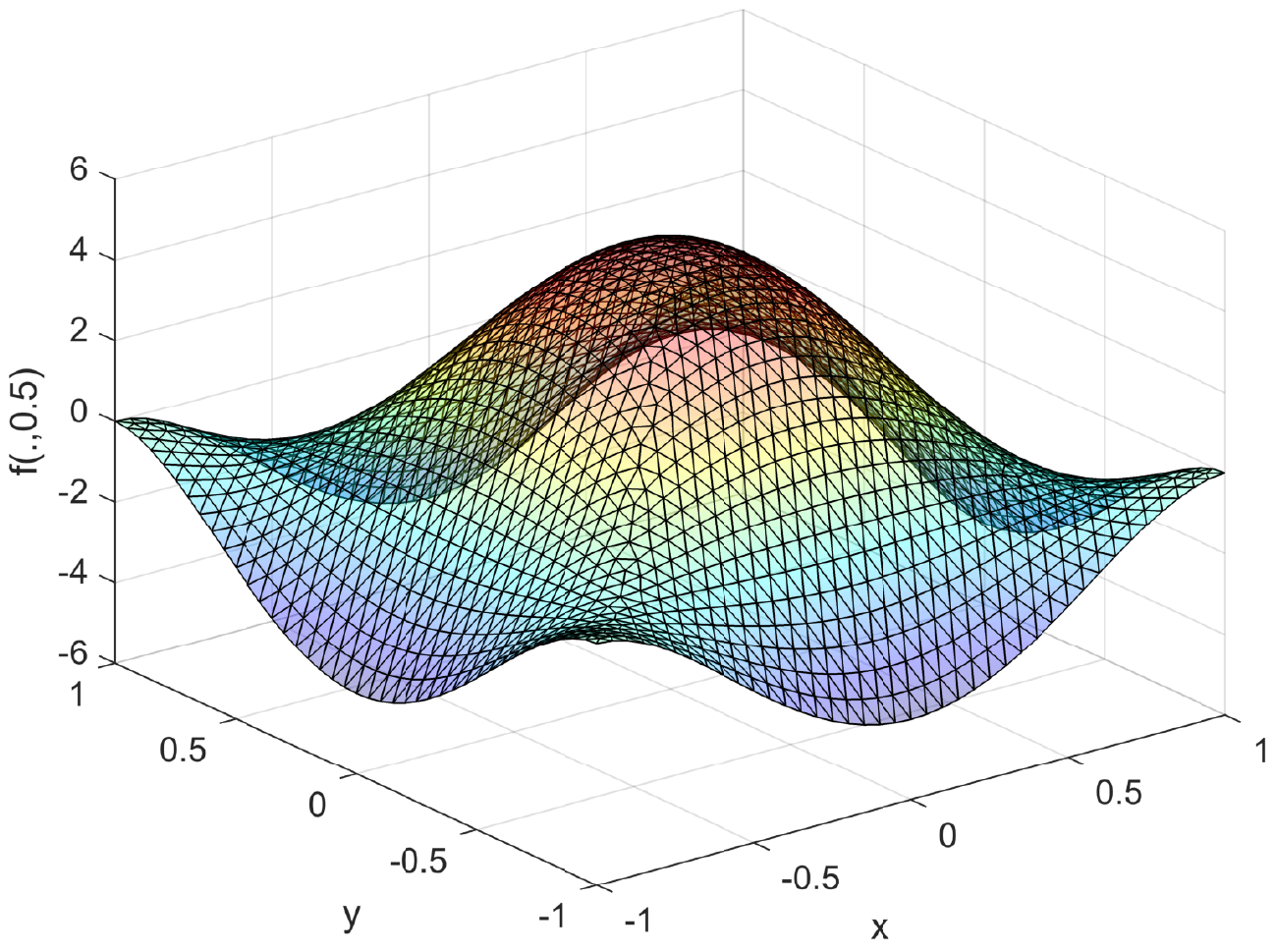}\\
			\hspace*{.45\textwidth} {\footnotesize (a)}
		\end{minipage}
		\begin{minipage}{.33\textwidth}
			\includegraphics[width=\textwidth]{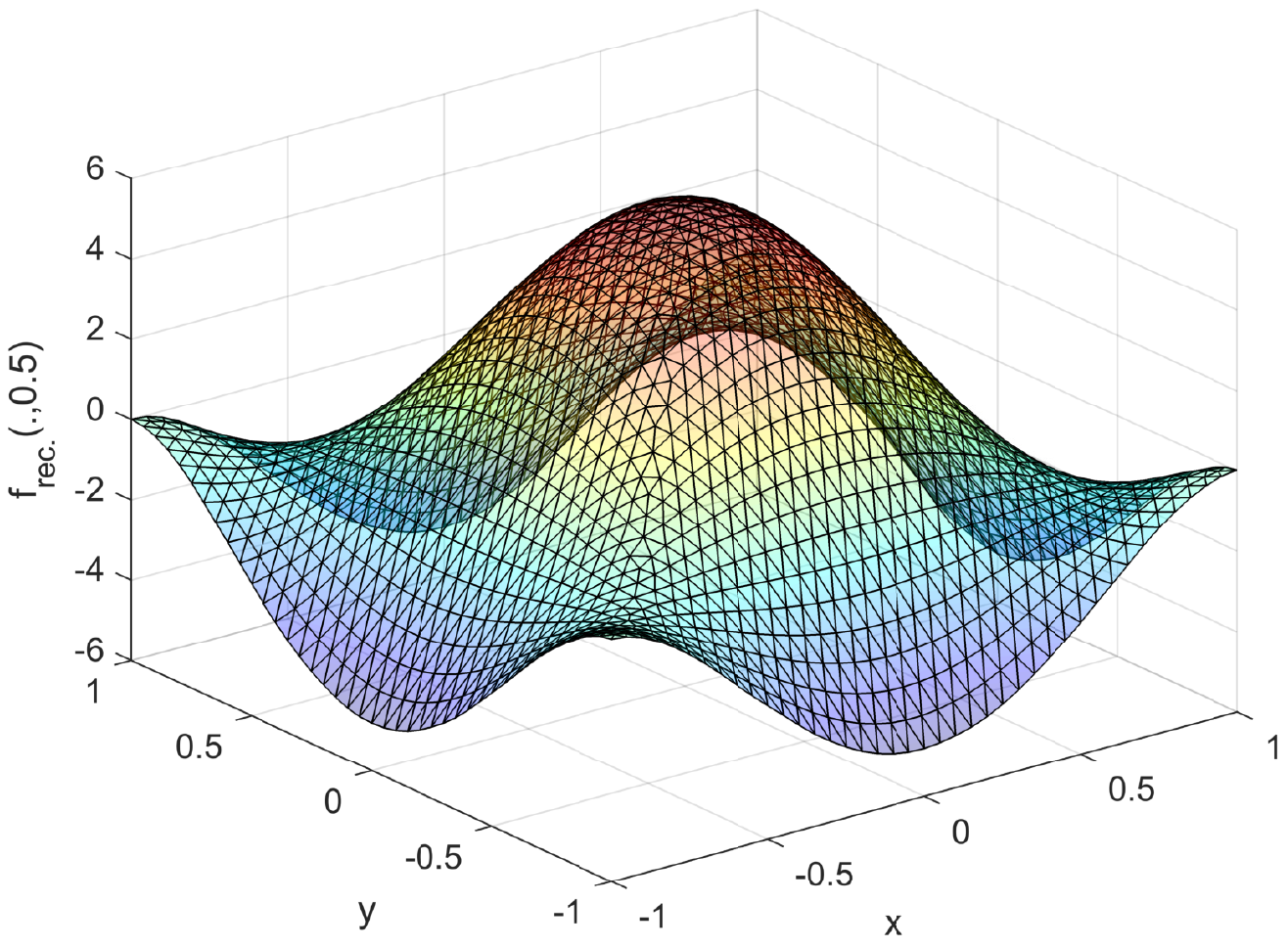} \\
			\hspace*{.45\textwidth} {\footnotesize (b)}
		\end{minipage}
	\begin{minipage}{.33\textwidth}
		\includegraphics[width=\textwidth]{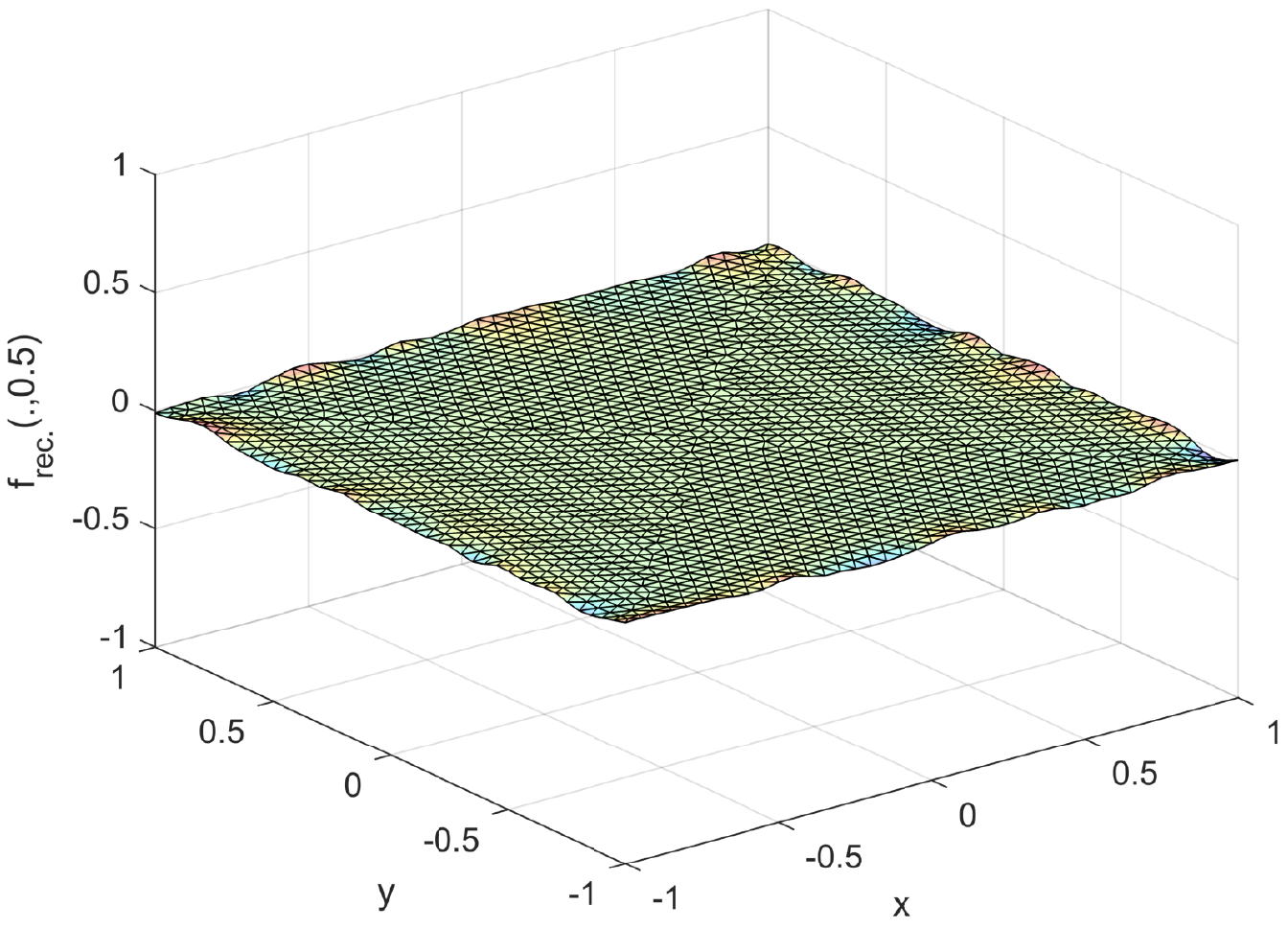}
		\hspace*{.45\textwidth} {\footnotesize (c)}
	\end{minipage}
	\end{center}
	\vspace*{-0.3cm}
	\caption{General source at $t = 0.5$: the exact source $f$ (a), recovered source (b), and recovered source with respect to the predicted source $f^*=0$.}
	\label{Fig:Mod1 comp_t05}
\end{figure}

In the example, we consider a general case, where $f$ depends on both space and time variables but with a somewhat simplified setting. 
We set $A = I_2$ the $2\times2$ identity matrix, $b = 0, \sigma = 0$ and $\Sigma = \partial\Omega \times (0,1]$. Moreover, we choose the exact state as $u(x,y,t) = t(x^2-1)^2(y^2-1)^2$. With this setting, it is straightforward to verify that the initial condition $q$, boundary condition $g$ and exact data $Z$ are all zeros. The exact source is then
\begin{equation*}
f_1(x,y,t) = (x^2-1)^2(y^2-1)^2 -t(x^2-1)^2(12y^2-4) -t(12x^2-4)(y^2-1)^2.
\end{equation*}
We in \Cref{Fig:Mod1 comp_P1} compare $f_1$ and $f_{rec.}$ at nodal point $P_1$, while in
\Cref{Fig:Mod1 comp_t05}(a) and \Cref{Fig:Mod1 comp_t05}(b) respectively perform  them at $t=0.5$. The numerical result shows that our reconstruction method produces a good approximation of the sought source in the general case.

We would like to discuss the role of the predicted source $f^*$ in the numerical solution to the identification problem. In principle, the CG method converges to the $f^*$-minimum-norm solution $f^\dag$. Therefore, for those inverse problems having more than one solution, the obtaining numerical solution approximates the one that is nearest to $f^*$. For the case considered here, we can easily observe that it accepts $f_2 = 0$ (corresponding the state $u_2=0$) as another solution. With the setting $f^*=0.5$ to imitate the situation where we have no information about the sought source, the algorithm converges to an approximation of $f_2$ as shown in \Cref{Fig:Mod1 comp_t05}(c).

\subsection{Identified source satisfying the condition \cref{17-6-20ct1}}
The last example aims at illustrating the convergence rate of the regularized approximations to the identification given in Theorem~\ref{con.rate}. To do so, the setting is similar to the previous subsections with the observations taking on the whole boundary instead. For generating the exact source, we start with the constant function $w=0.2$ on $\Sigma$ and then solve \eqref{10-12-19ct1} for the solution $F(w)$. As $w$ and the coefficients given in \cref{17-6-20ct2} are constant functions, it deduces that $F(w) \in C^\infty(\overline{\Omega} \times [0,T])$ (cf.\ \cite{evan,wolka}).  Let $f^* = (x^2+y)t$, we then have $f^\dagger = F(w)+f^* \in C^\infty(\overline{\Omega} \times [0,T])$. 

\begin{figure}[H]
	\begin{center}
		\begin{tabular}{ccc}
			\includegraphics[width=0.25\textwidth]{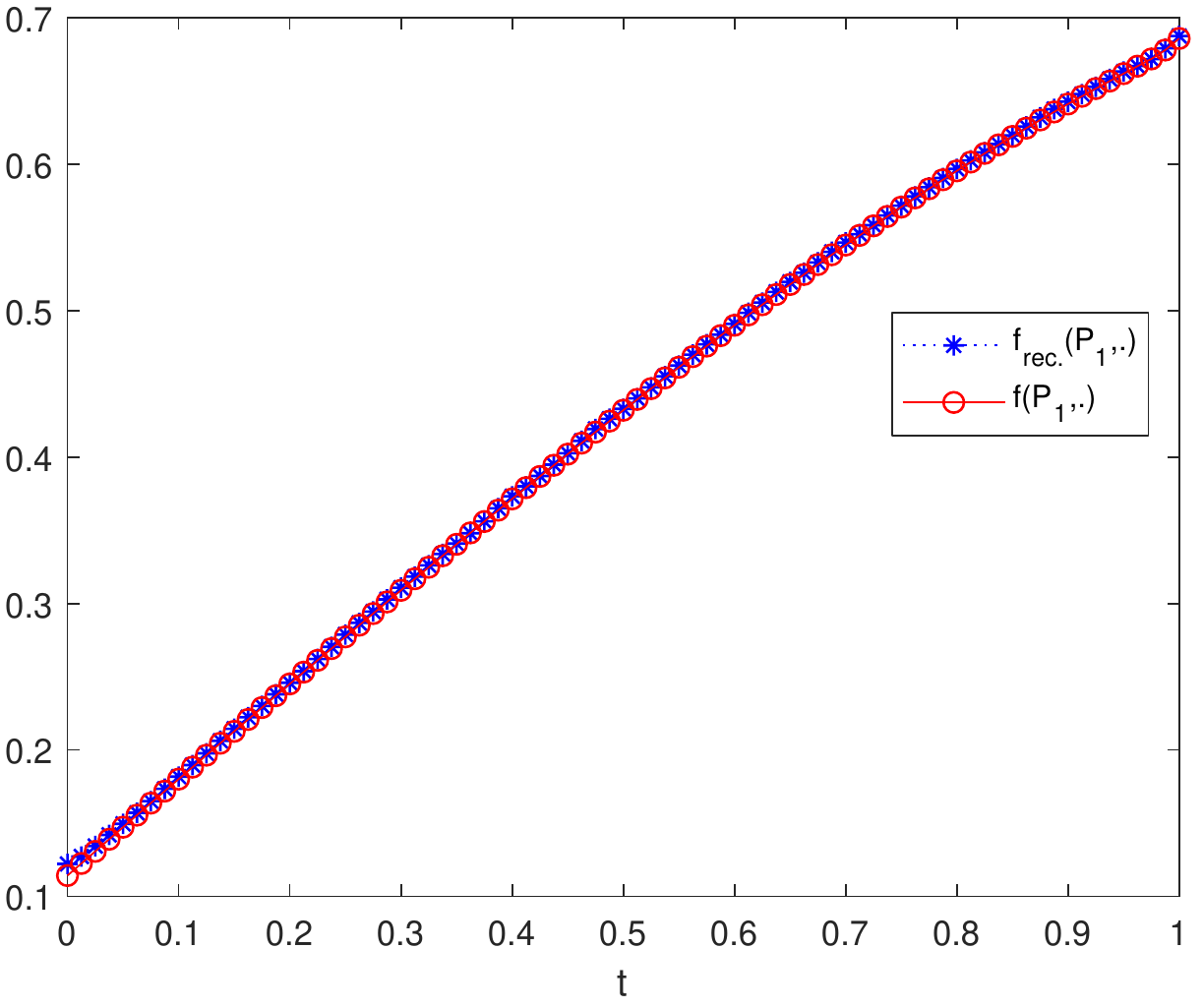} &
			\includegraphics[width=0.25\textwidth]{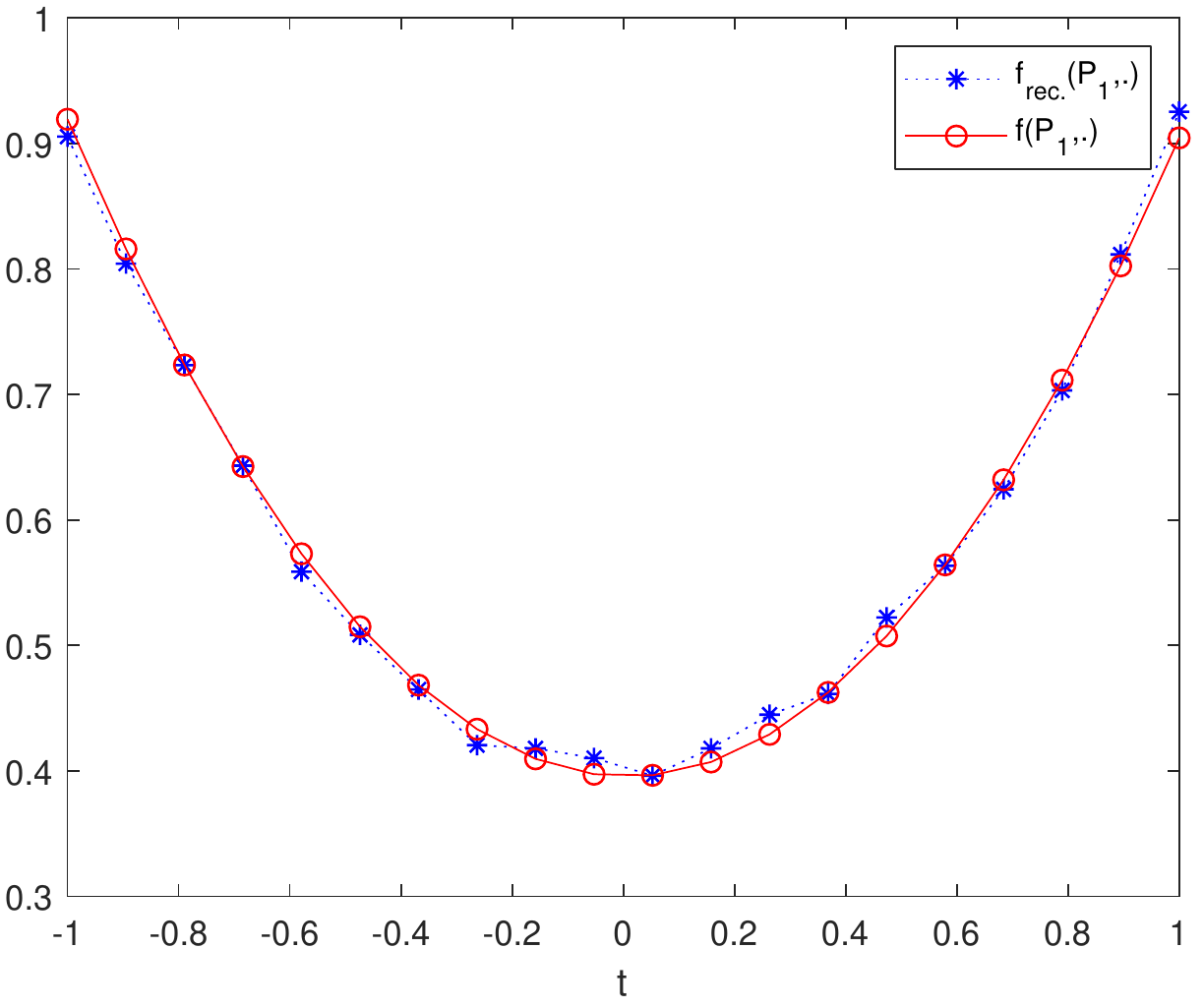} &
			\includegraphics[width=0.25\textwidth]{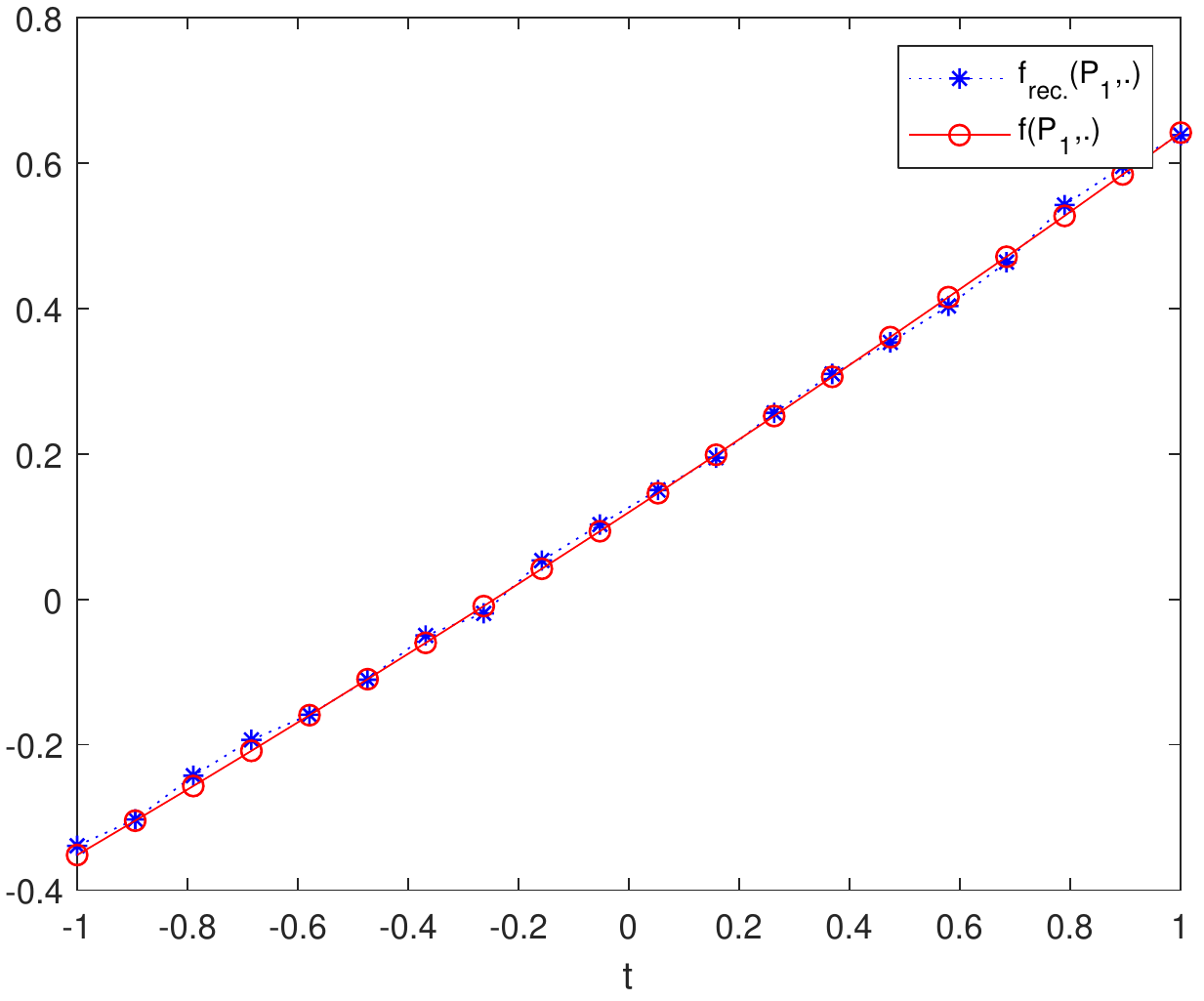}\\
			{\footnotesize (a)} & {\footnotesize (b)}& {\footnotesize (c)}
		\end{tabular}
		\begin{tabular}{cc}
			\includegraphics[width=0.35\textwidth]{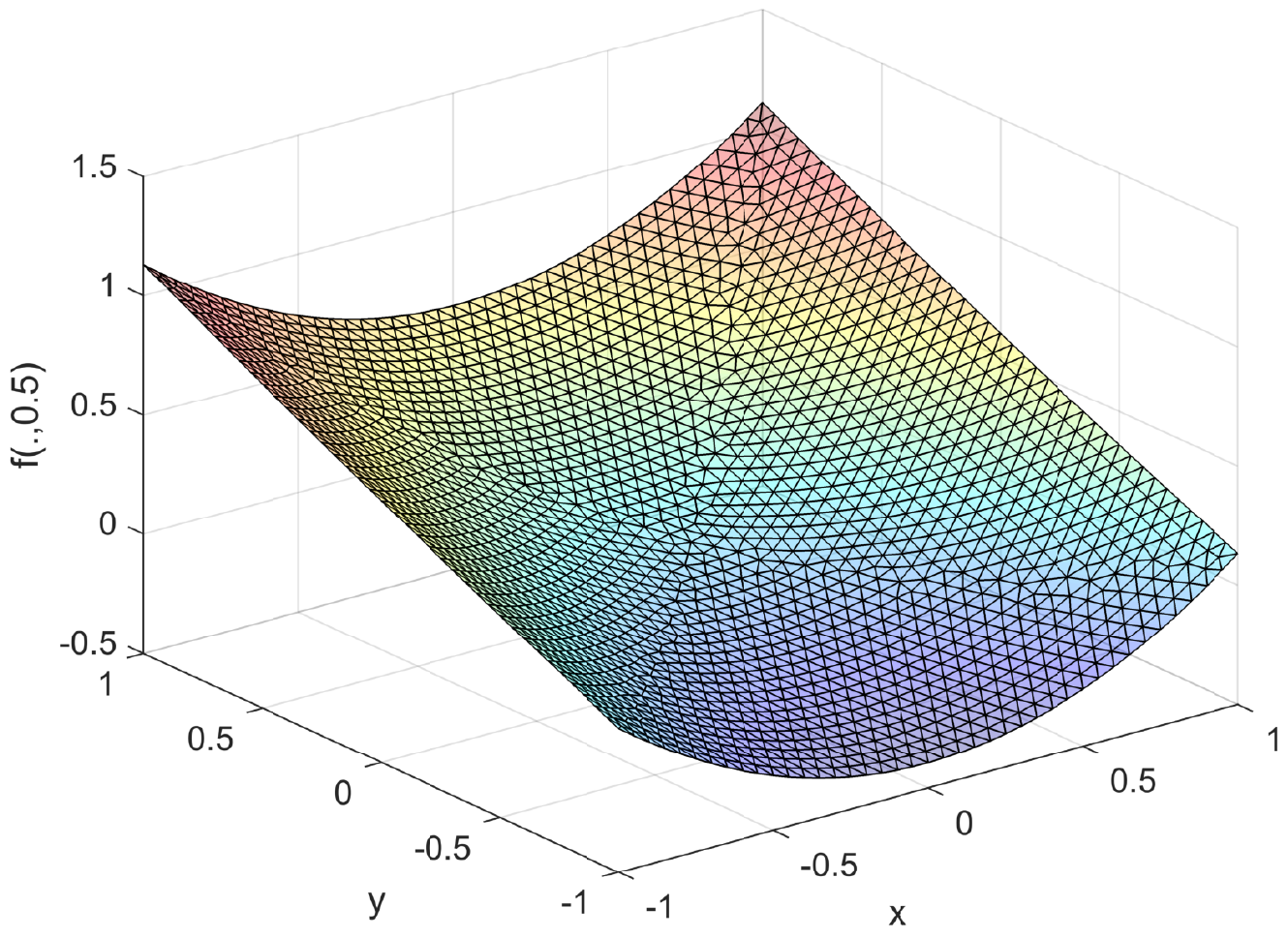} &
			\includegraphics[width=0.35\textwidth]{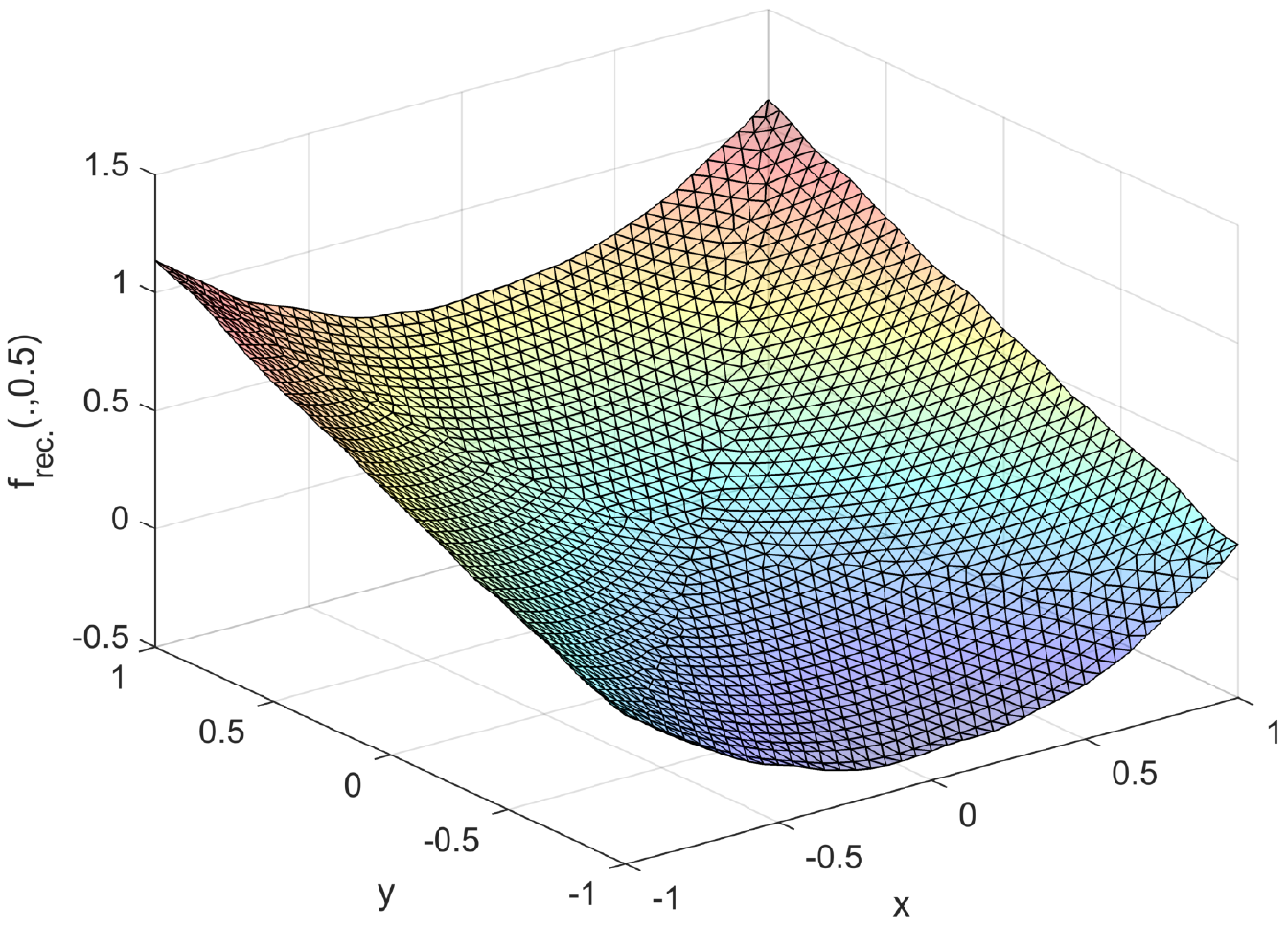}\\
			 {\footnotesize (d)} & {\footnotesize (e)}
		\end{tabular}
	\end{center}
	\caption{Identified source satisfying the condition \cref{17-6-20ct1}: comparisons along $t$ (a), $x$ (b), and $y$ (c) at $P_1$; the exact source $f$ (d), recovered source (e) at $t=0.5$.}
	\label{Fig:Mod4}
\end{figure}

In Table \ref{Tab:Mod4MeshRefine} and Table \ref{Tab:Mod4EOC}, we show the errors corresponding to different hierarchical mesh sizes and the resulting EOC, respectively, where all errors and
noisy levels get together smaller.
Figure \ref{Fig:Mod4} shows the comparison between the exact source and recovered one which match each other so well, as expected from our convergence result.

\begin{table}[H]
\caption{Identified source satisfying the condition \cref{17-6-20ct1}: Refinement level $l$, mesh size $h$, measurement error $\delta$, regularization parameter $\rho$, and errors}
	\begin{center}
\scalebox{0.95}{
	\begin{tabular}{ c l l l l l l}\hline
		$l$ & $h$ & $\delta$ & $\rho$& $\|u-u_{rec.}\|_{L^2(\Omega_T)}$& $\|u-u_{rec.}\|_{L^2(\Sigma)}$ &$\|f-f_{rec.}\|_{L^2(\Omega_T)}$\\  \hline
		1 & 0.8 & 0.32& 0.008&0.2199& 0.2709 &1.0994\\
		2 & 0.4 & 0.08& 0.004&0.0543& 0.0699 &0.3038\\
		3 & 0.2 & 0.02& 0.002&0.0136& 0.0175 &0.0842\\
		4 & 0.1 & 0.005&0.001&0.0032& 0.0042 &0.0237\\
		5 & 0.05 & 0.00125&0.0005&0.000745& 0.000929 &0.009616\\ \hline
	\end{tabular}}
	\end{center}
	\label{Tab:Mod4MeshRefine}
\end{table}

\begin{table}[H]
	\caption{Identified source satisfying the condition \cref{17-6-20ct1}: EOC}
	\begin{center}
	\begin{tabular}{ c l l l}\hline
		$l$ &  $\|u-u_{rec.}\|_{L^2(\Omega_T)}$ &$\|u-u_{rec.}\|_{L^2(\Sigma)}$ &$\|f-f_{rec.}\|_{L^2(\Omega_T)}$\\  \hline
		1 & -- & -- & -- \\
		2 & 2.0178& 1.9544 &1.8555 \\
		3 & 1.9973 &1.9979 & 1.8512\\
		4 & 2.0875&2.0589 & 1.8289\\
		5 & 2.1028&2.1766 & 1.3014\\
		Mean of EOC &2.0513 &2.0470& 1.7093 \\ \hline
	\end{tabular}
	\end{center}
	\label{Tab:Mod4EOC}
\end{table}

To close this section we wish to discuss about the Gibbs phenomenon that possibly appears when
numerically recovering of discontinuous functions concerning in this section, say Figures 1(c), 2(c), and 3(d). We observe that Gibbs phenomenon slightly happens in Figures 1(c) and 3(d), but it seems to be not in 2(c). These facts might be explained as follows. First, the use of differentiable regularization terms, e.g., the quadratic stabilizing penalty term, may smoothen the recovered solutions. We mention that to
reconstruct such discontinuous functions one usually employs the total variation regularization, which was originally introduced in image denoising \cite{ROF92}.
And second, by setting, in Figure 2(c) the a priori estimate $f^* =0$, while in the other cases, it is discontinuous. Also, obtaining numerical solutions performing in Figure 2(b) seem to be smooth, where the a priori estimate is differentiable; meanwhile we in 1(b) utilize a non-differentiable one.

\section{Conclusion}\label{appendix}

In this paper we investigate the inverse problem of identifying the source function in 
the parabolic equation
\begin{equation*}
\begin{aligned}
&\frac{\partial u}{\partial t} (x,t) + \mathcal{L}u(x,t) =f(x,t) \text{~in~} \Omega_T := \Omega \times (0,T],  \\
&\frac{\partial u(x,t)}{\partial\vec{n}} +\sigma(x,t)u(x,t) = g(x,t) \text{~on~} \mathcal{S} := \partial\Omega \times (0,T],\\
&u(x,0) = q(x) \text{~in~} \Omega
\end{aligned}
\end{equation*}
from a partial boundary measurement $z_\delta \in L^2(\Sigma)$ of the solution $u(x,t)$ on the surface $\Sigma := \Gamma\times(0,T) \subset\mathcal{S}$,
where $\mathcal{L}$ is a time-dependent, second order self-adjoint elliptic operator, $\Gamma$ is a relatively open subset of $\partial\Omega$ and $\delta>0$ is the error level of the observation.

The Crank-Nicolson Galerkin method is employed to fully discretize the parabolic equation.
As a result, the state $u(f)$ is then approximated by the finite sequence $(U^n_{h,\tau}(f))_{n=0}^M$ in which for each $n\in I_0$ the element 
$U^n_{h,\tau}(f) \in \mathcal{V}_h^1$ --- the space of piecewise linear,
continuous finite elements, where $h$ and $\tau$ is respectively the mesh size of the space discretization and the time step.

The least squares method combining with the quadratic stabilizing penalty term is utilized to tackle the identification problem, we then consider the unique minimizer $f_{\rho,\delta,h,\tau}$ of the minimization problem
$$
\min_{f\in L^2(\Omega_T)} \sum_{n=1}^M \int_{t^{n-1}}^{t^n} \|U^n_{h,\tau}(f)-z_\delta\|^2_{L^2(\Gamma)}dt + \rho\|f-f^*\|^2_{L^2(\Omega_T)} \eqno \left(\mathcal{P}_{\rho,\delta,h,\tau}\right)
$$
as a reconstruction, where $\rho>0$ is the regularization parameter and $f^*$ is an a priori estimate of the identified source.

We show that with $\delta,h,\tau$ approaching zero and an appropriate a priori regularization parameter choice $\rho=\rho(\delta,h,\tau)$
the whole sequence $\big(f_{\rho,\delta,h,\tau}\big)_{\rho>0}$ converges in the $L^2(\Omega_T)$-norm to the unique $f^*$-minimum-norm solution $f^\dag$ of the
identification problem as $\rho$ tends to zero.
The corresponding state sequence then converges in the $L^2(0,T; H^1(\Omega))$-norm to the state $u(f^\dag)$. Furthermore, the convergence rate \begin{align*}
\|f_{\rho,\delta,h,\tau} - f^\dag\|_{L^2(\Omega_T)} = \mathcal{O}(\sqrt{\delta})
\end{align*}
is established for an additional suitable source condition and an appropriate choice of the parameters $h$, $\tau$ and $\rho$ coupling with $\delta$. The numerical experiments are
presented to illustrate the efficiency of the theoretical findings.

\section*{Acknowledgments}
The authors would like to thank the Referees and the Editor for their valuable comments and
suggestions which helped to improve our paper.

T.N.T. Quyen gratefully acknowledges support of the University of Goettingen, Germany. N.T. Son is supported in part by
the Vietnam National Foundation for Science and Technology Development (NAFOSTED), grant 101.01-2017.319. The paper is completed when he is with the Universit\'e catholique de Louvain under the support of EOS Project no. 30468160 funded by FNRS and FWO.

\bibliographystyle{siam}

\end{document}